\documentclass[10pt]{amsart}

\usepackage{amssymb}
\usepackage{amsmath}
\usepackage{amsfonts}
\usepackage{geometry}
\usepackage{bbm}
\usepackage{hyperref}
\usepackage{tikz}
\usepackage{mathrsfs}
\usepackage[toc]{appendix}
\usepackage{multirow}
\usetikzlibrary{matrix,arrows,decorations.pathmorphing}
\usepackage{verbatim }
\usepackage{mathtools}
\usepackage{todonotes}
\usepackage{enumerate}

\newcounter{cprop}[section]

\newtheorem{theorem}[cprop]{Theorem}

\theoremstyle{plain}

\newtheorem{corollary}[cprop]{Corollary}

\newtheorem{lemma}[cprop]{Lemma}
\newtheorem{proposition}[cprop]{Proposition}

\newtheorem{assumption}[cprop]{Assumption}

\numberwithin{equation}{section}

\theoremstyle{definition}
\newtheorem{definition}[cprop]{Definition}

\theoremstyle{remark}
\newtheorem{remark}[cprop]{Remark}
\newtheorem{notation}[cprop]{Notation}

\renewcommand{\P}{\mathbb{P}}

\newcommand{\R}{\mathbb{R}}
\newcommand{\N}{\mathbb{N}}
\newcommand{\Z}{\mathbb{Z}}

\renewcommand{\d}{\mathrm{d}}

\newcommand{\vertiii}[1]{{\left\vert\kern-0.25ex\left\vert\kern-0.25ex\left\vert #1 
		\right\vert\kern-0.25ex\right\vert\kern-0.25ex\right\vert}}

\title[Inv. manifolds for random parabolic evolution equations with non-dense domain]{Invariant Manifolds for Random Parabolic Evolution Equations with almost sectorial operators}

\author{M. Ghani Varzaneh}
\address{Mazyar Ghani Varzaneh \\
	Fakult\"at für Mathematik und Informatik, FernUniversit\"at Hagen, Germany}
\email{mazyarghani69@gmail.com}

\author{F. Z. Lahbiri}
\address{Fatima Zahra Lahbiri \\
	Fakult\"at für Mathematik und Informatik, FernUniversit\"at Hagen, Germany}
\email{fatima-zahra.lahbiri@fernuni-hagen.de}

\author{S. Riedel}
\address{Sebastian Riedel \\
	Fakult\"at für Mathematik und Informatik, FernUniversit\"at Hagen, Germany}
\email{sebastian.riedel@fernuni-hagen.de}
\begin{document}
	
	\begin{abstract}
		In this paper, we develop a way of analyzing the random dynamics of stochastic evolution equations with a non-dense domain. Such problems cover several types of evolution equations. We are particularly interested in evolution equations with non-homogeneous boundary conditions of white noise type. We prove the existence of stable, unstable, and center manifolds around a stationary trajectory by combining integrated semigroup theory and invariant manifold theory. The results are applied to stochastic parabolic equations with white noise at the boundary.
		
		\vspace{0.2in}\noindent \textbf{Key words} Integrated semigroup, stochastic evolution equation, random dynamical systems, random invariant manifolds.
		
		\vspace{0.1in}\noindent \textbf{MSC:}  {\bf 47D06, 37H05, 47D62, 60H40}.
	\end{abstract}
	
	\maketitle 
	
	\section{Introduction}
	In this paper, we study dynamical properties of solutions to %It\^{o}'s
	stochastic evolution equations of the form
	\begin{align}\label{evolution-pb}
		\begin{cases}
			\d X(t)=AX(t) \, \mathrm{d}t+G(X(t)) \, \mathrm{d}t +  \mathrm{d}W(t), \quad t\geq0,\cr X(0)=\xi\in\overline{D(A)},
		\end{cases}
	\end{align}
	where $A:D(A)\subset H\to H$ is a linear operator on a separable Hilbert space $H$, $H_0:=\overline{D(A)}\neq H$, $G:H_0\to H$ is assumed to be globally Lipschitz continuous and 
	${W}$ represents a $Q$-Wiener process, with the covariance operator $Q \in \mathcal{L}(H)$.
	%defined on a complete filtered probability space $(\Omega,\mathcal{F},\mathbb{F},\mathbb{P})$, where $\mathbb{F}:=(\mathcal{F}_t)_{t\geq0}$ is the natural filtration. 
	Note that the Cauchy problem \eqref{evolution-pb} is well posed when $A$ is a densely defined Hille-Yosida operator. However, Cauchy problems with non-dense domains encompass an important range of differential equations, particularly those featuring nonhomogeneous boundary conditions. In fact, these boundary conditions are incorporated into supplementary state variables, which, consequently, make the evolution equation incapable of having a densely-defined linear operator. As a result, the well-known Hille-Yosida theory for $C_0$-semigroups is no longer applicable. In this article, our main focus lies on parabolic systems with non-homogeneous, random boundary data, i.e. white noise. For this case, recent studies have shown that well-posedness of the equation can still be obtained, see for instance \cite{Daprato-93,Sch-Ver11} and \cite{Ducrot-Lahbiri}.\smallskip
	% In this article, we study existence of random invariant manifolds for parabolic systems with non-homogeneous boundary data (white noise). 
	
	In the absence of noise, operators with non-dense domain have been studied using the concept of \emph{integrated semigroups}, which was first introduced by Arendt \cite{Arendt1,Arendt2} and further developed by Magal and Ruan \cite{Magal-Ruan} and Thieme \cite{Thieme1}. As mentioned before, parabolic systems with non-homogeneous boundary data can be reformulated using product spaces and matrix operators into the form of problem \eqref{evolution-pb} with a non-dense domain. In our case, the operator $A$ is subject to even weaker conditions than sectoriality. In fact, as pointed out in \cite{Prevost}, it is natural to consider an \emph{almost} sectorial operator $A$ (see Definition \ref{sectorial-operators}). \smallskip
	
	% The current work aims to study dynamical properties for parabolic systems with non-homogeneous boundary data (white noise). 
	% Using the integrated semigroup approach, we propose another take to the approach introduced in \cite{Daprato-Ergodicity}\todo{Put in the RDS section. Which class of equations did da Prato study?}. In fact the authors in \cite{Daprato-Ergodicity} used the dissipativity method to find conditions for existence and uniqueness of invariant measures and to establish the asymptotic properties of the corresponding transition semigroup. 
	% In this article, we study existence of random invariant manifolds for parabolic systems with non-homogeneous boundary data (white noise). 
	% Such problems are reformulated using product spaces and matrix operators into the form of problem \eqref{evolution-pb} with a non-dense domain, i.e., $\overline{D(A)} \neq H$. Furthermore, in this case, the operator $A$ is subject to even weaker conditions than sectoriality. In fact, as pointed out in \cite{Prevost} in the context of parabolic problems with non-homogeneous boundary conditions, it is natural to consider an \emph{almost} sectorial operator $A$ (see Definition \ref{sectorial-operators}).
	
	Our main goal in this work is to deduce the existence of random local invariant manifolds for solutions to \eqref{evolution-pb}. Note that center manifolds for non-densely evolution equations were investigated by Ruan and Magal in \cite{Magal-Mem09} in the deterministic setting. The authors used the Lyapunov-Perron method to study the existence and the smoothness of center manifolds for semilinear Cauchy problems with non-dense domain, which allows them to establish a Hopf bifurcation result for an age structure model. With a similar approach, Neam\c{t}u \cite{Neamtu2020} established the existence of stable and unstable manifolds for SPDE's with non-dense domain using integrated semigroups. Results regarding the existence of center manifold for SPDE's with non-dense domain can be found in the work of Li \& Zeng \cite{Li-Zeng22} where the authors used a similar method. \smallskip
	
	The main results in this work formulate conditions under which we can deduce the existence invariant manifolds (stable, unstable, and center) for a class of PDE's with white noise at the boundary, see Theorem \ref{stable_manifold}, Theorem \ref{unstable_manifold}, and Theorem \ref{center_manifold}. To prove these statements, we combine the multiplicative ergodic theorem (MET; see Theorem \ref{MET}) with the integrated semigroup approach. The MET was first proved in the literature by Oseledets \cite{Ose68} for matrix cocycles. Afterwards, it was generalized in various ways, see e.g, \cite{GVRS22,Blu16,GTQ15,MN02,Thi87,Rue82,Rag79,Rue79}. Note that the results obtained here are not only the first which provide the existence invariant manifolds for PDE's with white noise at the boundary, it's also the first attempt that combines the integrated semigroup approach with the multiplicative ergodic theorem.\smallskip
	% For compact cocycles the authors in \cite{GVR21} gave a proof based on measurable fields Banach spaces. Note that the results obtained here are not only the first which provides invariant manifolds for PDE's with white noise at the boundary, it's also the first attempt that combine integrated semigroup approach with random dynamical systems and multiplicative ergodic theorem.\\
	
	To formulate our main results, we use the theory of random dynamical systems (RDS) developed by L.~Arnold \cite{Arn98}. RDS form a powerful tool to study the long time behaviour of solutions to stochastic differential equations, both in finite and infinite dimensions. In fact, the objects that are described within RDS allow for a much finer analysis of the behaviour of stochastic processes compared to studying merely the invariant measure. This becomes crucial when studying possible bifurcations of stochastic differential equations. Random invariant manifolds are central objects in the theory of RDS. They are used to study e.g. stability, chaos, or stochastic bifurcation of stochastic differential equations. \smallskip
	
	Let us describe the strategy to prove our main theorems. The most critical property for a stochastic evolution equation to generate a random dynamical system is the cocycle property. A necessary condition for this to hold is the existence of a set of full measure on which the evolution equation induces a semiflow. For stochastic differential equations, the challenge lies in the fact that stochastic integrals notoriously generate nullsets that depend on the whole equation and, in particular, on the initial condition. In case the evolution equation is defined on a finite dimensional space (as it is the case for stochastic ordinary differential equations), the standard approach is to use Kolmogorov's continuity theorem to deduce the existence of a universal set of full measure on which the flow exists. If the state space is infinite dimensional, this approach breaks down. Our first main purpose in this work is to prove that the solution operator of the stochastic evolution equation \eqref{evolution-pb} generates a continuous random dynamical system, see Lemma \ref{dynamical-system} and Corollary \ref{dynamical-sys2}. This is done by a suitable decomposition of the SPDE into a linear equation for which the cocycle property can be deduced more easily and a random PDE that does not contain stochastic integrals. Then, under certain conditions, we prove that the cocycle is Fréchet differentiable and that linearizing the cocyle around a stationary point yields a linear cocycle acting in the closed state space $H_0$, see Proposition \ref{Fréchet}. For the linearized cocycle, we have to check that the integrability condition of the multiplicative ergodic theorem holds. The initial step towards achieving this involves establishing a suitable a priori bound (see Theorem \ref{Priory-bound-cocycle}, Corollary \ref{estimate on differnce} and Proposition \ref{PPR}). We can then apply abstract results in the literature (in particular \cite{GVR23} and \cite{Sebastian-maz23}) that allow us to deduce the existence of random local stable, unstable and center manifolds.  
	%     Note that to get the stationary solution, we need to make sure the function $G$ stays within certain limits. It should either be globally bounded (see Lemma \ref{Example-G-bounded}) everywhere or meet a specific condition (like estimate \ref{condition on G}, see Lemma \ref{integrability_derivative}). This is important in analyzing the long-time behavior of our problem. If $G$ behaves wildly or grows without control, the solution won't work well. So, we have to pay close attention to how $G$ behaves to make sure our solution stays on track.
	% Wconsider $X$  a Banach and we denote by $\mathcal{L}\left( X,H\right) $
	% 	the space of bounded linear operators from $X$ into $H$ and by $%
	% 	\mathcal{L}\left( H\right) $ the space $\mathcal{L}\left( H,H\right) .$
	
	\smallskip

	\section{Preliminaries}
	In this section, we first collect some basic results about integrated semigroups from \cite{Ducrot-Magal-Prevost,Magal-Ruan}. Let ${A}$ be a linear operator on a Banach space ${X}$ with a non-dense domain, its resolvent set is denoted by $\rho({A}):=\{\lambda\in\mathbb{C}:\lambda I-{A}\quad\text{is invertible} \}$ and its spectrum is denoted by $\sigma({A}):=\mathbb{C}\backslash\rho({A})$. Note that $\lambda\in\rho({A})$ if $R(\lambda,{A}):=(\lambda I-{A})^{-1}\in\mathcal{L}({X})$, where by $\mathcal{L}(X)$, we mean the space of linear operators from $X$ to $X$, equipped with the usual operator norm. Now, let ${X}_{0}=\overline{D({A})}$ and denote the part of ${A}$ on ${X}_{0}$ by ${A}_{0}:D({A}_{0})\subset {X}_{0}\to {X}_{0}$, which is a linear operator defined by
	\begin{align}\label{Def-domaine}
		{A}_{0}x={A}x,\quad D({A}_{0}):=\{y\in D({A}):{A}y\in X_{0}\}.
	\end{align}
	Assume that there exists a constants $\omega$ satisfying $(\omega,+\infty)\subset\rho({A})$, then for each $\lambda>\omega$,
	\begin{equation*}
		D({A}_{0})=R(\lambda,{A}){X}_{0},\quad R(\lambda,{A}_{0})=R(\lambda,{A})|_{{X}_{0}}.
	\end{equation*}
	It follows from \cite[Lemma 2.1]{Magal-Ruan} that $\rho({A})=\rho({A}_{0})$, and thus $\sigma({A})=\sigma({A}_{0})$. The following result from \cite[Lemma 2.2.11]{Magal-Ruan} gives necessary and sufficient conditions to ensure that
	\begin{equation*}
		\overline{D({A}_{0})}=\overline{D({A})}.
	\end{equation*}
	\begin{lemma}
		Let ${A}:D({A})\subset{X}\to{X}$ be a linear operator on $({X},\|\cdot\|)$. Assume that there exist two constants, $\omega\in\mathbb{R}$ and $M>0$, such that $(\omega,+\infty)\subset\rho({A})$ and for each $\lambda>\omega$
		\begin{equation*}
			\limsup_{\lambda\to+\infty}\lambda\left\|R(\lambda,{A}_{0})\right\|=\limsup_{\lambda\to+\infty}\lambda\left\|R(\lambda,{A})\right\|_{\mathcal{L}({X}_{0})}<+\infty.
		\end{equation*}
		Then the following properties are equivalent:
		\begin{itemize}
			\item [$(i)$] $\lim_{\lambda\to +\infty}\lambda R(\lambda,{A})x=x,$ $\forall x\in{X}_{0}$;
			\item [$(ii)$] $\lim_{\lambda\to +\infty} R(\lambda,{A})x=0,$ $\forall x\in{X}$;
			\item [$(iii)$] $\overline{D({A}_{0})}={X}_{0}$.
		\end{itemize}
	\end{lemma}
	The definition of an integrated semigroup is given as follows.
	\begin{definition}
		Let $({X},\|\cdot\|)$ be a Hilbert space. A family of bounded linear operator $({S}(t))_{t\geq0}$ is called an \textbf{integrated semigroup} on ${X}$ if the following properties are satisfied:
		\begin{itemize}
			\item [$(i)$] ${S}(0)=0$;
			\item [$(ii)$] $t\to{S}(t)x$ is continuous from $[0,+\infty)$ into ${X}$ for each $x\in{X}$;
			\item [$(iii)$] ${S}(t)$ satisfies
			\begin{equation*}
				{S}(t){S}(s)=\int_{0}^{t}\left[{S}(s+r)-{S}(r)\right] \, \d r,\quad\forall t,s\in [0,+\infty).
			\end{equation*}
		\end{itemize}
		An integrated semigroup $({S}(t))_{t\geq0}$ on ${X}$ is said to be {\bf nondegenerate} if
		\begin{equation*}
			{S}(t)x=0,\forall t\geq0\Rightarrow x=0.
		\end{equation*}
		An integrated semigroup $({S}(t))_{t\geq0}$ on ${X}$ is said to be {\bf exponentially bounded} if there exist two constants $\tilde{M}\geq1$ and $\tilde{\omega}>0$, such that
		\begin{equation*}
			\left\|{S}(t)\right\|_{\mathcal{L}({X})}\leq\tilde{M}e^{\tilde{\omega}t},\quad\forall t\geq0.
		\end{equation*}
	\end{definition}
	
	In the following, we will suppose that the operator ${A}$ satisfies the conditions of the Hille-Yosida theorem with the exception of the density of $D({A})$ in ${X}$. Note that ${A}$ is said to be a Hille-Yosida operator on ${X}$ if there exist two constants $M\geq 1$ and $\omega\in\R$  such that $(\omega,+\infty)\subset\rho({A})$ and
	\begin{equation*}
		\left\|R(\lambda,{A})^{k}\right\|_{\mathcal{L}({X})}\leq\frac{M}{(\lambda-\omega)^{k}},\quad\forall\lambda>\omega,\forall k\geq 1.
	\end{equation*}
	Concerning ${A}$, we make the following assumptions.
	\begin{assumption}
		\label{Assumption-A0} We assume that
		\begin{itemize}
			\item [$(i)$] ${A}$ is a Hille-Yosida operator on ${X}_{0}$;
			\item [$(ii)$] $\lim_{\lambda\to +\infty} R(\lambda,{A})x=0,$ $\forall x\in{X}$.
		\end{itemize}
	\end{assumption}
	Under assumption \ref{Assumption-A0} and according to \cite[Lemma 3.4.2]{Magal-Ruan}, the operator ${A}_{0}$ generates a strongly continuous semigroup (called also $\mathcal{C}_0$-semigroup) on ${X}_{0}$ denoted by $({T}(t))_{t\geq0}$. Furthermore, the operator ${A}$ generates an integrated semigroup on ${X}$ denoted by $({S}(t))_{t\geq0}$ (see \cite[Proposition 3.43]{Magal-Ruan}). Now we are ready to state some important properties of integrated semigroups.
	\begin{proposition}
		Let assumption \ref{Assumption-A0} be satisfied. Then ${A}$ generates a non-degenerate integrated semigroup $({S}(t))_{t\geq0}$ that satisfies for each $x\in{X}$, each $t\geq0$, and each $\nu>\omega$
		\begin{align}\label{Forme-integrated-sg}
			\begin{split}
				{S}(t)&=(\nu I-{A}_{0})\int_{0}^{t}{T}(s) \, \d s R(\nu,{A}),\\
				{S}(t)x&=\nu\int_{0}^{t}{T}(s)R(\nu,{A})x\, \d s + \left[I-{T}(t)\right]R(\nu,A)x.
			\end{split}
		\end{align}
		In addition, the map $t\to{S}(t)$ is continuously differentiable if and only if $x\in{X}_{0}$ and
		\begin{equation*}
			\frac{d{S}(t)x}{dt}={T}(t)x,\quad\forall t\geq0,\forall x\in{X}_{0}.
		\end{equation*}
	\end{proposition}
	Next, we recall the basic existence and uniqueness result of the non-homogenous Cauchy problem in \cite{Magal-Ruan} for
	\begin{equation}\label{non-homogenous-CP}
		\frac{\d z}{\d t} = {A}z(t)+f(t),\quad t\in[0,\tau],\quad z(0)=x\in{X}_{0}.
	\end{equation}
	\begin{proposition}
		For each $f\in C^{1}(0,\tau,{X})$,  set
		\begin{equation*}
			({S}*f)(t):=\int_{0}^{t}{S}(s)f(t-s)\,  \d s,\quad\forall t\in [0,\tau].
		\end{equation*}
		Then the following properties hold:
		\begin{itemize}
			\item [$(i)$] The map $t\to({S}*f)(t)$ is continuously differentiable on $[0,\tau]$;
			\item [$(ii)$] $({S}*f)(t)\in D({A}),\quad\forall t\in[0,\tau]$;
			\item [$(iii)$] For each $\lambda>\omega$, and $t\in[0,\tau]$, we have
			\begin{equation*}
				R(\lambda,{A})\frac{d}{dt}({S}*f)(t)=\int_{0}^{t}{T}(t-s)R(\lambda,{A})f(s)\, \d s.
			\end{equation*}
		\end{itemize}
	\end{proposition}
	Besides, as concluded in \cite{Magal-Ruan}, the non-homogenous Cauchy problem \eqref{non-homogenous-CP} has a unique mild solution $z\in C([0,\tau];{X}_{0})$ given by
	\begin{equation*}
		z(t)={T}(t)x + \lim_{\lambda\to+\infty}\int_{0}^{t}{T}(t-s)\lambda R(\lambda,{A})f(s) \, \d s,\quad\forall t\in[0,\tau].
	\end{equation*}
	Let us recall the notion of sectorial operators.
	
	\begin{definition}
		The operator ${A}_{0}$ is said to be {\bf sectorial} if there are constants $\omega\in\mathbb{R}$, $\theta\in ]\frac{\pi}{2},\pi[$, and $\bar{M}>0$ such that
		\begin{itemize}
			\item [$(i)$] $\rho({A}_{0})\supset \mathcal{S}_{\theta}=\{z\in\mathbb{C}; z\neq 0|\arg (z-\bar{\omega})|\leq\theta\}$;
			\item [$(ii)$] $\|R(\lambda,{A}_{0})\|\leq\frac{\bar{M}}{|\lambda-\bar{\omega}|},\quad\lambda\in\mathcal{S}_{\theta}.$
		\end{itemize}
	\end{definition}
	Note that when $A_0$ is a sectorial operator, then $A_0$ is the infinitesimal generator of a strongly continuous analytic semigroup $({T}(t))_{t\geq0}$ on $X_0$. According to \cite[Theorem 6.8]{Pazy-book}, we can define for each $\beta\geq0$ the operator $(-A_0)^{-\beta}:D((-A_0)^{-\beta})\subset X_0\to X_0$ as the inverse of $(-A_0)^{\beta}$. Moreover, we know (see \cite[Theorem 6.13]{Pazy-book}) that for each $t>0$, $(-A_0)^{-\beta}T(t)$ is a bounded operator and 
	\begin{equation}\label{crucial-estimate1}
		\|(-A_0)^{-\beta}T(t)\|\leq M_{\beta}t^{-\beta}e^{\omega_A t}.
	\end{equation}
	
	We finally introduce the notion of almost sectorial operators, which is a weaker notion than sectoriality
	
	\begin{definition}\label{sectorial-operators}
		The operator ${A}$ is said to be {\bf $\beta-$almost sectorial} if there are constants $\omega\in\mathbb{R}$, $\theta\in ]\frac{\pi}{2},\pi[$, and $\bar{M}>0$ such that
		\begin{itemize}
			\item [$(i)$] $\rho({A})\supset \mathcal{S}_{\theta}=\{z\in\mathbb{C}; z\neq 0|\arg (z-\bar{\omega})|\leq\theta\}$;
			\item [$(ii)$] $\|R(\lambda,{A})\|\leq\frac{\bar{M}}{|\lambda-\bar{\omega}|^{\beta}},\quad\lambda\in\mathcal{S}_{\theta}.$
		\end{itemize}
	\end{definition}
	For parabolic equations, the existence and uniqueness of solutions for \eqref{non-homogenous-CP} has been studied in \cite{Ducrot-Magal-Prevost} when the operator ${A}$ is $\frac{1}{p^*}$-almost sectorial, ${A}_{0}$ is sectorial for each $x\in {X}_{0}$ and $f\in L^{p^*}([0,\tau],{X})$. Moreover, it has been proved (see \cite[Lemma 3.7 and Lemma 3.8]{Ducrot-Magal-Prevost}) that
	the map $t\to{S}(t)$ is continuously differentiable and
	\begin{equation}\label{Crucial-equlity-isg}
		\frac{dS(t)x}{dt}=(-A_0)^{\beta}T(t)(-A)^{-\beta}x,\quad\forall t>0,\quad\forall x\in X,\quad\forall\beta>1-\frac{1}{p^*}
	\end{equation}
	where 
	\begin{align}\label{fractional_growth}
		\Vert (-A)^{-\beta}\Vert<\infty,\quad\forall\beta>1-\frac{1}{p^*}.
	\end{align}
	Moreover, for $q^*\in [1,+\frac{1}{1-\frac{1}{p^*}})$ and for $\tau>0$,
	\begin{align}\label{Crucial-inequlity-isg}
		\int_{0}^{\tau}\left\|\frac{dS(t)}{dt}\right\|^{q^*}dt\leq (M_\beta \|((-A)^{-\beta})\|)^{q^*}\int_{0}^{\tau}t^{-q^*}e^{-\omega_A t}dt<\infty,
	\end{align}
	see also, \cite[Lemma 9.1.8]{Magal-Ruan}.\smallskip

	Next, we introduce some basic notions from the theory of random dynamical systems. For more details, the reader is referred to \cite{Arn98}. Let $(\Omega,\mathcal{F},\P)$ be a probability space and let $\mathcal{T}$ be either $\Z$ or $\R$. Let $\theta_t \colon \Omega \to \Omega$, $t \in \mathcal{T}$, be a family of measurable maps. We will assume that the family $(\theta_t)_{t \in \mathcal{T}}$ is $\mathbb{P}$-preserving, i.e $\theta_t\mathbb{P}=\mathbb{P}$ for all $t \in \mathcal{T}$. Furthermore, we assume that
	\begin{enumerate}[(i)]
		\item The mapping $(t,\omega)\to\theta_t \omega$ is $(\mathcal{B}(\mathcal{T}) \otimes \mathcal{F}, \mathcal{F})$ measurable;
		\item $\theta_{0}\omega=\omega$;
		\item $\theta_{t+s}=\theta_{t}\circ\theta_{s}$ for all $t,s\in\mathcal{T}$.  
	\end{enumerate}
	The quadruple $(\Omega,\mathcal{F},\mathbb{P},(\theta_t)_{t\in\mathcal{T}})$ is called {\bf metric dynamical system}. If $\mathcal{T} = \Z$, we set $\theta \coloneqq \theta_1$. If $\mathbb{P}$ is ergodic with respect to $(\theta_t)_{t\in\mathcal{T}}$, i.e. any invariant subset has zero or full measure, we call the system {\bf ergodic}.\smallskip
	
	Set $\mathcal{T}_+ \coloneqq \{t \in \mathcal{T}\, |\, t \geq 0\}$. A \textbf{cocycle} on $H_0$ is a map
	\begin{align*}
		\phi \colon \mathcal{T}_+ \times \Omega \times H_0 \to H_0 
	\end{align*}
	that is $(\mathcal{B}(\mathcal{T}_+) \otimes \mathcal{F} \otimes \mathcal{B}(H_0) ,\mathcal{B}(H_0))$-measurable and satisfies the {\bf cocycle property}, i.e. there exists a measurable set of full measure $\tilde{\Omega}$ such that
	\begin{enumerate}[(i)]
		\item $\phi(0,\omega,x) = x$ for every $x \in H_0$ and $\omega \in \tilde{\Omega}$ and
		\item $\phi(t+\tau,\omega,x)=\phi(t,\theta_{\tau}\omega,\phi(\tau,\omega,x))$ for every $t, \tau \in \mathcal{T}_+$, $\omega \in \tilde{\Omega}$ and $x \in H_0$.
	\end{enumerate}
	The cocycle property presents a generalization of the semigroup property. In particular, if the $\omega$-dependence is omitted, then the cocycle property reduces to the semigroup property. A metric dynamical system together with a cocycle is called \textbf{random dynamical system} (RDS). The cocycle $\phi$ is called \textbf{continuous} if the map
	\begin{align*}
		\phi(\cdot,\omega,\cdot) \colon \mathcal{T}_+ \times H_0 \to H_0
	\end{align*}
	is continuous for every $\omega \in \tilde{\Omega}$. It is called \textbf{linear} (resp. {\bf compact} or \textbf{(Fr\'echet-)differentiable}), if the map 
	\begin{align*}
		\phi(t,\omega,\cdot) \colon H_0 \to H_0
	\end{align*}
	is linear (resp. compact or (Fr\'echet-)differentiable) for every $\omega \in \tilde{\Omega}$ and every $t \in \mathcal{T}_+$. If the cocycle map is continuous (resp. linear, compact or (Fr\'echet-)differentiable), the corresponding random dynamical system will be called continuous (resp. linear, compact or (Fr\'echet-)differentiable), too.
	
	\section{Existence of the Random Dynamical System and Well-Posedness}
	In this section, we show that \eqref{evolution-pb} generates a random dynamical system by a suitable decomposition of the equation. Let us explain the idea. We first consider the linear problem
	\begin{align*}
		\d Y(t) &= AY(t) \, \d t + \d W(t), \quad t \geq 0, \\
		Y(0) &= Y_0 \in H_0.
	\end{align*}
	This equation can explicitly be solved for any inital condition, even for random $Y_0$. Then we solve the random PDE
	\begin{align*}
		\d V(t) &= AV(t) \, \d t + G(V(t) + Y(t))\, \d t, \quad t \geq 0, \\
		V(0) &= \xi - Y_0 \in H_0.
	\end{align*}
	An informal calculation reveals that the solution $X(t)$ to \eqref{evolution-pb} can be given by $X(t) = V(t) + Y(t)$. Indeed, $X(0) = V(0) + Y(0) = \xi - Y_0 + Y_0 = \xi$ and
	\begin{align*}
		\d X(t) &= \d V(t) + \d Y(t) \\
		&= AV(t) \, \d t + G(V(t) + Y(t))\, \d t + AY(t) \, \d t + \d W(t) \\
		&= A(V(t) + Y(t)) \, \d t + G(V(t) + Y(t))\, \d t + \d W(t) \\
		&= AX(t) \, \d t + G(X(t)) \, \d t + \d W(t).
	\end{align*}

	We assume that $\{(W(t)\}_{t \in \R}$ is a two-sided $Q$-Wiener process. Let's first consider the linear part of \eqref{evolution-pb} (i.e $G=0$)
	\begin{align}\label{Linear-Problem}
		\begin{cases}
			\partial_t Y(t)=A Y(t)+\mathrm{d}W(t),\quad t\in[0,\tau],\\
			Y(0)=\xi\in H_0.
		\end{cases}
	\end{align}
	The main assumption can be simplified as follows:
	\begin{assumption}
		\label{Assumption1.1} We assume that
		\begin{itemize}
			\item[(a)] the operator $A_{0}$, the part of $A$ in $\overline{D(A)}$, is sectorial.
			\item[(b)] There exist $\omega_A \in \mathbb{R}$ and $p^{\ast }\in \left[
			1,+\infty \right)$ such that $\left( \omega_A ,+\infty \right) \subset \rho
			\left( A\right) $, the resolvent set of $A$, and
			\begin{equation*}
				\underset{\lambda \rightarrow +\infty }{\limsup }\;\lambda ^{\frac{1}{%
						p^{\ast }}}\left\Vert \left( \lambda I-A\right) ^{-1}\right\Vert _{\mathcal{L%
					}\left( H\right) }<+\infty. 
			\end{equation*}  
		\end{itemize}
	\end{assumption}
	This assumption is equivalent to $\rho(A)\neq\emptyset$. Since $\rho(A)$ and $\rho(A_0)$ are non-empty we have $\rho(A)=\rho(A_0)$ (see \cite[Lemma 2.1]{Magal-Mem09}). Note that the assumption \ref{Assumption1.1}-$(a)$ ensures that the operator $A_0$ is a sectorial Hille-Yosida operator on $\overline{D(A)}$ (see \cite[Lemma 2.1 and Lemma 2.2]{Magal-Ruan}). In this case, $A_{0}$, the part of $A$ in $\overline{D(A)}$ generates a strongly continuous analytic semigroup  $\left\{ T_{A_{0}}(t)\right\} _{t\geq 0}$ on $\overline{D(A)}$. According \cite[Lemma 3.4.2, Proposition 3.4.3]{Magal-Ruan} the operator $A$ generates an analytic integrated semigroup on $H$ denoted by $({S}_A(t))_{t\geq0}$ (see ) such that for all $\tau>0$
	\begin{equation}\label{regularity-integrated-sg}
		t\mapsto S_A(t)\in W^{1,q^*}(0,\tau;\mathcal L(X))\text{ with }q^*\in \left[1,\frac{1}{1-\frac{1}{p^*}}\right).
	\end{equation}
	%			\item[(c)] As the family $({S}_A(t))_{t\geq0}$ of linear operators turns out to be differentiable with respect to $t>0$, then we assume that $\left\{\frac{\mathrm{d}S_A(\cdot)}{\mathrm{d}t}\right\}$ is a Hilbert-Schmidt operator from $H$ to itself.
	It has been proven that under the above assumptions, the mild solution (or integrated solution) of \eqref{Linear-Problem} exists and is given by 
	\begin{align*}
		Y(t) &\coloneqq {T}_{0}(t)\xi+\frac{\d}{\mathrm{d}t}\int_{0}^{t}{S}_A(t-s)\, \d {W}(s),\\
		&= {T}_{0}(t)\xi + \lim_{\lambda\to +\infty} \int_{0}^{t}{T}_{0}(t-s)\lambda R(\lambda,{A})\, \d {W}(s),
	\end{align*}
	for any $t\in [0,\tau]$ and $\mathbb{P}$-a.s. Consider the following Cauchy problem with non-dense domain
	\begin{align}\label{Problem-2}
		\begin{cases}
			\partial_t Y_\omega(t) = A Y_\omega (t) + \mathrm{d}W(t),\quad t\in[0,\tau],\\
			Y_\omega(0) =Y_\omega \in H_0
		\end{cases}
	\end{align}
	where
	\begin{equation*}
		Y_\omega = \int_{-\infty}^{0}\frac{\d}{\d \tau}S_A(-\tau) \, \mathrm{d} W(\tau).
	\end{equation*}
	If $\omega_{A}<0$, then for any $t\in [0,\tau]$ the integrated solution of the Cauchy problem \eqref{Problem-2} can be rewritten as
	\begin{align}\label{eqn:def_Y}
		Y_{\omega}(t) = \int_{-\infty}^{t}\frac{\d}{\d \tau} S_A(t-\tau) \, \mathrm{d}W(\tau),\quad \mathbb{P}-a.s.
	\end{align}
	Note that according to \eqref{regularity-integrated-sg} and as $Q$ has finite trace, the stochastic integral \eqref{eqn:def_Y} is well defined.
	\begin{assumption}
		\label{Assumption1.2}  We assume that $\theta:\mathbb{R}\times\Omega\to\Omega$ is a family of $\mathbb{P}-$preserving transformations where $(\Omega,\mathcal{F},\mathbb{P},\{\theta_t\}_{t\in\mathbb{R}})$ is an ergodic metric dynamical system 
		%           %and $\theta_t\omega=\omega(\cdot+t)-\omega(t)$ is the Wiener-shift 
		that satisfies
		\begin{equation*}
			W_t(\theta_s\omega)=W_{t+s}(\omega)-W_{s}(\omega)
		\end{equation*}
		for every $s,t \in \R$ and every $\omega \in \Omega$.
		% \begin{itemize}
			% 	\item[(a)] $\theta:\mathbb{R}\times\Omega\to\Omega$ is a family of $\mathbb{P}-$preserving transformations where $(\Omega,\mathcal{F},\mathbb{P},\theta_{t\in\mathbb{R}})$ is a metric dynamical system 
			%           %and $\theta_t\omega=\omega(\cdot+t)-\omega(t)$ is the Wiener-shift 
			%           that satisfies
			% 	\begin{equation*}
				% 		W_t(\theta_s\omega)=W_{t+s}(\omega)-W_{s}(\omega).
				% 	\end{equation*}
			% 	\item[(b)] $Y_\omega$ is stationary, i.e
			% 	\begin{equation}
				% 		Y_{\theta_t\omega}(s)=Y_\omega(t+s).  \label{1.2}
				% 	\end{equation}
			% \end{itemize}
	\end{assumption}
	It is well-known that one can always construct a $Q$-Wiener process that satisfies Assumption \ref{Assumption1.2}.
	
	\begin{lemma}\label{lemma:Y_stationary}
		Let assumptions \ref{Assumption1.1} and \ref{Assumption1.2} hold. Then the random process $Y_\omega(t)$ defined in \eqref{eqn:def_Y} is \emph{stationary}, i.e. there is a measurable set of full measure $\tilde{\Omega}$ such that for every $\omega \in \tilde{\Omega}$ and $s,t \in \R$,
		\begin{equation}
			Y_{\theta_t\omega}(s)=Y_\omega(t+s).  \label{1.2}
		\end{equation}
		% 	\begin{equation*}
			% 		\mathbb{E}\| Y_{\theta_t\omega}(0)\|^{2}<\infty.
			% 	\end{equation*}
	\end{lemma}
	
	\begin{proof}
		For $s$ and $t$ fixed, we have
		\begin{align*}
			Y_{\theta_t\omega}(s) &= \int_{-\infty}^{s}\frac{\d}{\d \tau} S_A(s-\tau) \, \mathrm{d}W_{\tau}(\theta_t\omega) \\
			&= \int_{-\infty}^{s}\frac{\d}{\d \tau} S_A(s-\tau) \, \mathrm{d}W_{\tau + t}(\omega) \\
			&=  \int_{-\infty}^{t + s}\frac{\d}{\d \tau} S_A(t + s-\tau) \, \mathrm{d}W_{\tau}(\omega) \\
			&= Y_{\omega}(t + s)
		\end{align*}
		almost surely, i.e. on a set of full measure $\Omega_{s,t}$ that may depend on $s$ and $t$. Our goal is to find a universal set of full measure on which this identity holds. First, since the integrated semigroup $(S_A(t))_{t\geq0}$ satisfies \eqref{regularity-integrated-sg}, it follows from \cite{Da-Za} that the process $t\to Y_\omega (t)$ is continuous. According to Kolmogorov's theorem the two processes $(s,t)\to Y_{\theta_\omega t}(s)$ and $(s,t)\to Y_{\omega}(t+s)$ also have continuous modifications.
		%Initially, the equality \ref{1.2} holds only a set of a full measure $\Omega_{s,t}$ that depends on $s$ and $t$.
		% , however we need this equality to hold for all $\omega$ and for all combinations of $s$ and $t$. To do so, let's $(s,t)$ be a pair of rational numbers. As the set of all pairs of rational numbers is countable, 
		We can consider the intersection of all events $$\Omega^{\prime} \coloneqq \bigcap_{(s,t)\in\mathbb{Q}^2}\Omega_{s,t}$$ where $\mathbb{P}(\Omega^{\prime})=1$, i.e, the equality \ref{1.2} holds for every $\omega\in\Omega^{\prime}$ and any pair of rational numbers $(s,t)$. Using the density of rational numbers, for any real numbers $(s,t)\in\mathbb{R}$, there exist sequences of rational numbers ${s_n}$ and ${t_n}$ converging to $s$ and $t$ respectively. Finally, the continuity properties of the stochastic processes allow us to take the limit, ensuring that \ref{1.2} holds for all real $s$ and $t$ on $\Omega^{\prime}$.
	\end{proof}
	% \begin{remark}
		% 	A sufficient condition for the random variable $Y_\omega$ to be stationary is
		% 	\begin{equation*}
			% 		\mathbb{E}\| Y_{\theta_t\omega}(0)\|^{2}<\infty.
			% 	\end{equation*}
		% \end{remark}

	\begin{proposition}\label{result1}
		Under assumptions \eqref{Assumption1.1}-\eqref{1.2}, the integrated solution corresponding to the linear Cauchy problem \eqref{Linear-Problem} generates a random dynamical system $\phi:\mathbb{R^+}\times\Omega\times H_0\to H_0$.
	\end{proposition}
	
	\begin{proof}
		Let $\omega\in\Omega$, $t>0$ and define the system $\phi:\mathbb{R^+}\times\Omega\times H_0\to H_0$ by
		\begin{equation}\label{dynamical-system}
			\phi_{\omega}^{t}(\xi):=T_{0}(t)\xi-\frac{\d}{\mathrm{d}t}S_{A}(t)Y_{\omega}(0)+Y_{\theta_{t}\omega}(0),
		\end{equation}
		where $Y_{\omega}(\cdot)$ is the solution of \eqref{Problem-2}. It suffices to prove that $\phi$ satisfies the cocycle property, i.e,
		\begin{equation*}
			\phi_{\omega}^{t+s}(\xi)= \phi_{\theta_{t}\omega}^{s}\circ\phi_{\omega}^{t}(\xi), ~~\text{for all}~~ \xi\in H_0,~t,s\in\mathbb{R}^{+}~~\text{and all}~~\omega\in\Omega.
		\end{equation*}
		The expression on the left-hand side can be straightforwardly expressed by
		\begin{equation*}
			\phi_{\omega}^{t+s}(\xi):=T_0(t+s)\xi-\frac{\d}{\mathrm{\d}t}S_{A}(t+s)Y_\omega+Y_{\theta_{t+s}\omega}(0).
		\end{equation*}
		As the integrated solution satisfies $Y_{\theta_{t}\omega}(0)\in H_0$, the right-hand side yields
		\begin{align*}
			\phi_{\theta_{t}\omega}^{s}\circ\phi_{\omega}^{t}(\xi)&= T_0(s)(\phi_{\omega}^{t}\xi)-\frac{\d}{\mathrm{d}s}S_{A}(s)Y_{\theta_{t}\omega}(0)+Y_{\theta_{s}\theta_{t}\omega}(0)\\
			&=T_0(s)(T_{0}(t)\xi-\frac{\d}{\mathrm{d}t}S_{A}(t)Y_{\omega}(0)+Y_{\theta_{t}\omega}(0))-\frac{\d}{\mathrm{d}s}S_{A}(s)Y_{\theta_{t}\omega}(0)+Y_{\theta_{t+s}\omega}(0)\\
			&=T_0(t+s)\xi-\frac{\d}{\mathrm{d}t}S_{A}(t+s)Y_\omega(0)+T_0(s)Y_{\theta_{t}\omega}(0)-\frac{\d}{\mathrm{d}s}S_{A}(s)Y_{\theta_{t}\omega}+Y_{\theta_{t+s}\omega}(0)\\
			&=T_0(t+s)\xi-\frac{\d}{\mathrm{d}t}S_{A}(t+s)Y_\omega+Y_{\theta_{t+s}\omega}(0)
		\end{align*}
		and the cocycle property holds.
	\end{proof}
	We now claim that the solution to \eqref{evolution-pb} can be written as $X(t) = V(t) + Y_{\theta_t \omega}(0)$ where $V$ solves the following random PDE
	\begin{align}\label{eqn:random_PDE}
		\begin{cases}
			\mathrm{d}V(t) = AV(t)\, \mathrm{d}t + G(V(t) +Y_{\theta_t \omega}) \, \mathrm{d}t, \quad t\in[0,\tau],\cr V(0)=\xi - Y_\omega\in H_0,
		\end{cases}
	\end{align}
	and $Y_\omega(t)$ is an integrated solution to \eqref{Problem-2} which belongs to $\mathcal{C}([0,\tau],H_0)$. To prove the well-posedness of \eqref{eqn:random_PDE} it suffices to prove that it has a fixed point. The result is as follows.
	\begin{theorem}\label{Existence-solution}
		Let assumptions \eqref{Assumption1.1} and \eqref{Assumption1.2} be satisfied with $\omega_A<0$. Then the equation \eqref{eqn:random_PDE} admits a unique integrated solution (or mild solution) $V(\cdot)\in\mathcal{C}([0,\tau],H_0)$ for every $\tau > 0$.
	\end{theorem}
	
	\begin{proof} Let $f\in L^{2}([0,\tau],H)$. The proof relies on a fixed point argument. Let
		\begin{equation*}
			V(t) \coloneqq T_{0}(t)(\xi-Y_\omega) +\int_{0}^{t}\frac{\mathrm{d}}{\mathrm{d}s}S_{A}(t-s)f(s) \, \mathrm{d}s.
		\end{equation*}
		According to \cite[Lemma 3.7]{Ducrot-Magal-Prevost}, for each $t > 0$, $\beta \in \left(\frac{1}{q^*}, 1\right)$, and $x \in H$
		\begin{equation*}
			\frac{\mathrm{d}}{\mathrm{d}t}S_{A}(t)x=(-A_0)^{\beta}T_{0}(t)(-A)^{-\beta}x.
		\end{equation*}
		Therefore, we have for $0\leq t_1\leq t_2<\tau$,
		\begin{align*}
			&\|V(t_2)-V(t_1)\|_{H_0}\\
			&= \bigg\|(T_0(t_2)-T_0(t_1))(\xi-Y(0))+\int_{t_1}^{t_2}(-A_0)^{\beta}T_{0}(t_2-s)(-A)^{-\beta}f(s) \, \mathrm{d}s\\
			&\quad +\int_{0}^{t_1}(T_0(t_2-t_1)-I)(-A_0)^{\beta}T_{0}(t_1-s)(-A)^{-\beta}f(s) \, \mathrm{d}s \bigg\|_{H_0}\\
			&\leq \|(T_0(t_2-t_1)-I)T_0(t_1)(\xi-Y(0))\|_{H_0}+\left\|\int_{t_1}^{t_2}(-A_0)^{\beta}T_{0}(t_2-s)(-A)^{-\beta}f(s) \, \mathrm{d}s\right\|_{H_0}\\
			&\quad +\left\|\int_{0}^{t_1}(T_0(t_2-t_1)-I)(-A_0)^{\beta}T_{0}(t_1-s)(-A)^{-\beta}f(s) \, \mathrm{d}s\right\|_{H_0}.
		\end{align*}
		Since
		\begin{equation*}
			\|(T_0(t_2-t_1)-I)T_0(t_1)(\xi-Y(0))\|\leq C \|\xi-Y(0)\|_{H_0},
		\end{equation*}
		we get that
		\begin{equation}\label{equality1}
			\lim_{t_2\to t_1}\|(T_0(t_2-t_1)-I)T_0(t_1)(\xi-Y(0))\|_{H_0}=0.
		\end{equation}
		Furthermore, by H\"older's inequality,
		\begin{align}\label{equality2}
			&\lim_{t_2\to t_1}\left\|\int_{t_1}^{t_2}(-A_0)^{\beta}T_{0}(t_2-s)(-A)^{-\beta}f(s) \, \mathrm{d}s\right\|^2_{H_0}\nonumber\\
			&\leq\lim_{t_2\to t_1}\int_{t_1}^{t_2}\|(-A_0)^{\beta}T_{0}(t_2-s)(-A)^{-\beta}\|^2_{\mathcal{L}(H)} \, \mathrm{d}s\int_{t_1}^{t_2}\|f(s)\|^2_{H} \, \mathrm{d}s\\
			&\leq C \lim_{t_2\to t_1}\int_{t_1}^{t_2}\|f(t)\|^2_{H} \, \mathrm{d}t=0. \nonumber
		\end{align}
		It is clear that
		\begin{equation*}
			\left\|\int_{0}^{t_1}(-A_0)^{\beta}T_{0}(t_1-s)(-A)^{-\beta}f(s) \, \mathrm{d}s\right\|_{H_0}^2\leq C \int_{0}^{t_1}\|f(s)\|^{2}_{H} \, \mathrm{d}s.
		\end{equation*}
		From Lebesgue’s dominated convergence theorem, it follows that
		\begin{align}\label{equality3}
			\lim_{t_2\to t_1}\left\|(T_0(t_2-t_1)-I)\int_{0}^{t_1}(-A_0)^{\beta}T_{0}(t_1-s)(-A)^{-\beta}f(s) \, \mathrm{d}s\right\|_{H_0}=0.
		\end{align}
		From \eqref{equality1}-\eqref{equality2} and \eqref{equality3}, we obtain that
		\begin{equation*}
			\lim_{t_2\to t_1}\|V(t_2)-V(t_1)\|_{H_0}=0.
		\end{equation*}
		Thus, we proved that $V(\cdot)\in\mathcal{C}([0,\tau],H_0)$.\\
		For fixed $\tau>0$, we define the map $\mathcal{J}: \mathcal{C}([0,\tau],H_0)\to \mathcal{C}([0,\tau],H_0)$ for any $V(\cdot)\in\mathcal{C}([0,\tau],H_0)$ by
		\begin{equation*}
			\mathcal{J}(V)(t) \coloneqq T_0(t)\xi + \int_{0}^{t}\frac{\d}{\mathrm{d}s}S_{A}(t-s) \, \mathrm{d}W(s) + \int_{0}^{t}\frac{\d}{\mathrm{d}s}S_{A}(t-s)G(V(s)+Y_{\theta_t}\omega) \, \mathrm{d}s.
		\end{equation*}
		Hence under assumption \eqref{Assumption1.1}, $\mathcal{J}$ is well defined. Then we can find a constant $c(\tau)>0$ such
		\begin{align*}
			|\mathcal{J}(V_1)-\mathcal{J}(V_2)|_{\mathcal{C}([0,\tau],H_0)}^2&\leq c(\tau)|G(V_1+Y_{\theta_\cdot}\omega)-G(V_2+Y_{\theta_\cdot}\omega)|_{L^{2}([0,\tau],H_0)}^2\\
			&\leq  c(\tau) L|V_1-V_2|_{L^{2}([0,\tau],H_0)}^2\\
			&\leq  c(\tau) L|V_1-V_2|_{\mathcal{C}([0,\tau],H_0)}^2
		\end{align*}
		where for any $t\in [0,\tau]$, $\lim_{t\to 0}c(t)=0$. Assume $\tau$ is small enough such that $c(\tau) L<1$, thus $\mathcal{J}$ is contractive if $c(\tau)<\frac{1}{L}$. By means of the Banach fixed point theorem, the random PDE \eqref{eqn:random_PDE} (and hence\eqref{evolution-pb}) has a unique integrated solution, which is continuous in $t$ for every $\omega\in\Omega$. For large $\tau>0$, we can now subdivide the interval $[0,\tau]$ into sufficiently small pieces on which the solutions exist and glue them together to obtain a global solution.
	\end{proof}
	
	\begin{corollary}\label{dynamical-sys2}
		Let assumptions \eqref{Assumption1.1} and \eqref{Assumption1.2} be satisfied and assume that $\omega_A<0$. Then the mild solution of \eqref{evolution-pb} generates a random dynamical system $\tilde{\phi}:\mathbb{R}^{+}\times\Omega\times H_0\to H_0$.
	\end{corollary}
	\begin{proof}
		For any $t\geq0$, $\omega\in\Omega$ and $\xi\in H_0$, we consider the system
		\begin{align}\label{cocycle-NP}
			\tilde{\phi}^t_{\omega}(\xi):=\psi^t_{\omega}(\xi-Y(\omega))+Y_{\theta_t \omega},
		\end{align}
		where
		\begin{equation*}
			\psi(0,t,\omega)(\xi-Y(\omega)):=\psi^t_{\omega}(\xi-Y(\omega))=T_{0}(t)(\xi-Y(\omega))+\int_{0}^{t}(-A_0)^{\beta}T_{0}(t-s)(-A)^{-\beta}G(\psi^s_{\omega}(\xi-Y(\omega))+Y_{\theta_s \omega}) \, \mathrm{d}s.
		\end{equation*}
		As the random PDE \eqref{eqn:random_PDE} satisfies the flow property, by using Proposition \ref{result1} it suffices to show that \eqref{eqn:random_PDE} generates the dynamical system $\psi:\mathbb{R^+}\times\Omega\times H_0$. A sufficient condition in order to so is to show that the following equality holds
		\begin{equation*}
			\psi(u,t,\omega)=\psi(0,t-u,\theta_{u}\omega),\quad 0<u<t.
		\end{equation*}
		Now, let $\xi \in H_0$ and $0<u<t$, we have
		\begin{align*}
			\psi(0,t-u,\theta_{u}\omega)(\xi-Y(\omega))&=T_{0}(t-u)(\xi-Y(\omega))\\&+\int_{0}^{t-u}(-A_0)^{\beta}T_{0}(t-u-\tau)(-A)^{-\beta}G(\psi(0,\tau,\theta_{u}\omega)(\xi-Y(\omega))+Y_{\theta_{\tau+u}\omega})d\tau.
		\end{align*}
		In the other hand
		\begin{align*}
			&\psi(u,t,\omega)(\xi-Y(\omega))=T_{0}(t-u)(\xi-Y(\omega))+\int_{u}^{t}(-A_0)^{\beta}T_{0}(t-s)(-A)^{-\beta}G(\psi(u,s,\omega)(\xi-Y(\omega))+Y_{\theta_{s}\omega})\mathrm{d}s.\\
			&=T_{0}(t-u)(\xi-Y(\omega))+\int_{0}^{t-u}(-A_0)^{\beta}T_{0}(t-u-\tau)(-A)^{-\beta}G(\psi(u,\tau+u,\omega)(\xi-Y(\omega))+Y_{\theta_{\tau+u}\omega})d\tau.
		\end{align*}
		We find that
		\begin{align*}
			\psi(0,t-u,\theta_{u}\omega)(\xi-Y(\omega))-\psi(u,t,\omega)(\xi-Y(\omega))&=\int_{0}^{t-u}(-A_0)^{\beta}T_{0}(t-u-\tau)(-A)^{-\beta}[G(\psi(0,\tau,\theta_{u}\omega)(\xi-Y(\omega))\\&+Y_{\theta_{\tau+u}\omega})-G(\psi(u,\tau+u,\omega)(\xi-Y(\omega))+Y_{\theta_{\tau+u}\omega})]d\tau
		\end{align*}
		Using \eqref{Crucial-inequlity-isg} together with Gronwall’s inequality, the desired equality follows.
	\end{proof}
	
	\section{Differentiability, compactness, and some a priori bounds}
	
	The key to proving the existence of invariant manifolds lies in establishing several a priori bounds. A crucial step will be to show that the linearized cocycle satisfies a mild integrability property. To see this, we will first prove that under natural assumptions on the coefficients, the solution to \eqref{evolution-pb} is indeed differentiable in the initial condition.
	
	\begin{proposition}\label{Fréchet}
		Let assumptions \ref{Assumption1.1} and \ref{1.2} be satisfied and assume that $\omega_A<0$. Let $G:H_0\to H$ be of class $\mathcal{C}^1$. Then the solution $\tilde{\phi}^t_{\omega}(\xi)$ of \eqref{evolution-pb} is Fréchet differentiable at $\xi$. Moreover, for any $\beta\in(1-\frac{1}{p^*},1)$ and $t>0$, the derivative at $\xi\in H_0$ in the direction $\eta\in H_0$ satisfies
		\begin{equation}\label{Fréchet_differentiable}
			D_{\xi}\tilde{\phi}^t_{\omega}[\eta]=T_{0}(t)\eta+\int_0^t (-A_0)^{\beta}T_{0}(t-s)(-A)^{-\beta}DG_{\tilde{\phi}^s_{\omega}(\xi)}(D_{\xi}\tilde{\phi}^t_{\omega}[\eta]) \, \mathrm{d}s.
		\end{equation}
	\end{proposition}
	
	\begin{proof}
		For any $t>0$ and $\xi\in H_0$, we set
		\begin{equation*}
			v(t):=\tilde{\phi}^t_{\omega}(\xi)-\xi
		\end{equation*}
		and define the map $\mathcal{F}:H_0\times\mathcal{C}([0,\tau_0],H_0)\to\mathcal{C}([0,\tau_0],H_0)$ by
		\begin{equation*}
			\mathcal{F}(h,v) \coloneqq T_{0}(t)(v(0)+h)+\int_0^t (-A_0)^{\beta}T_{0}(t-s)(-A)^{-\beta}G(v(s)+h) \, \mathrm{d}s+Y_{\theta_{t}\omega}.
		\end{equation*}
		Using similar reasoning as in Theorem \ref{Existence-solution}, we establish that for $h\in H_0$ and $v\in \mathcal{C}([0,\tau_0],H_0)$, the map $\mathcal{F}$ admits a fixed point $v$, i.e.,
		\begin{equation*}
			\mathcal{F}(h,v)=v.
		\end{equation*}
		This can be equivalently reformulated by defining the function $\tilde{\mathcal{F}}$ as
		\begin{equation*}
			\tilde{\mathcal{F}}(h,v):=\mathcal{F}(h,v)-v=0.
		\end{equation*}
		To prove that $\tilde{\phi}^t_{\omega}(\xi)$ is Fréchet differentiable at $\xi$, we will apply the Implicit Function Theorem. First, we need to verify the necessary conditions. Specifically, we must show that the derivative of $\tilde{\mathcal{F}}(h,v)$ with respect to $v$ is an isomorphism. A sufficient condition for this is to verify whether  
		\[
		\left\|\frac{\partial}{\partial v} \tilde{\mathcal{F}}(h,v)\right\| < 1.
		\]
		We proceed as follows. As $\mathcal{F}$ is a  contractive map, i.e., there exists a constant $K<1$ such that 
		\begin{equation*}
			\|\mathcal{F}(h,v)-\mathcal{F}(h,\tilde{v})\|\leq K\|v-\tilde{v}\|.
		\end{equation*}
		In addition, if we consider $\tilde{v}:=v+R$ we obtain
		\begin{equation*}
			\lim_{\|R\|\to 0}\frac{\|\mathcal{F}(h,v)-\mathcal{F}(h,\tilde{v})\|}{\|R\|}=\left\|\frac{\partial}{\partial v}\tilde{\mathcal{F}}(h,v)\right\|\leq K<1.
		\end{equation*}
		For all $\eta,\xi\in H_0$, we set
		\begin{align*}
			P^t(\epsilon,\eta,\omega):=\frac{\tilde{\phi}^t_{\omega}(\xi+\epsilon\eta)-\tilde{\phi}^t_{\omega}(\xi)}{\epsilon}.
		\end{align*}
		We find that
		\begin{align*}
			P^t(\epsilon,\eta,\omega)&=T_{0}(t)\eta+\frac{1}{\epsilon}\int_0^t (-A_0)^{\beta}T_{0}(t-s)(-A)^{-\beta}G(\tilde{\phi}^s_{\omega}(\xi+\epsilon\eta))-G(\tilde{\phi}^s_{\omega}(\xi)) \, \mathrm{d}s\\
			&=T_{0}(t)\eta+\frac{1}{\epsilon}\int_0^t (-A_0)^{\beta}T_{0}(t-s)(-A)^{-\beta}\int_0^1 DG_{A(\theta,\xi,\eta,s)}(\phi^s_{\omega}(\xi+\epsilon\eta)-\tilde{\phi}^s_{\omega}(\xi)) \, \d \theta \, \mathrm{d}s,
		\end{align*}
		where 
		\begin{equation*}
			A(\theta,\xi,\eta,s):=\tilde{\phi}^s_{\omega}(\xi)+\theta(\tilde{\phi}^s_{\omega}(\xi+\epsilon\eta)-\tilde{\phi}^s_{\omega}(\xi)).
		\end{equation*}
		Since $G$ is assumed to be of class $\mathcal{C}^1$ and a Lipschitz continuous function, then it has a bounded growth, i.e., there exists $r>0$ such that
		\begin{equation}\label{Derivative-of-G}
			|x|<r,\quad |D_{x}G|< k(r).
		\end{equation}
		Using Gronwall's inequality, we get for any $\eta\in H_0$ and $t>0$
		\begin{align*}
			\|P^t(\epsilon,\eta,\omega)\|&\leq Me^{\omega_A t}\|\eta\|_{H_0}+C_{\beta} \left(\int_0^t \|\int_0^1 DG_{A(\theta,\xi,\eta,s)}(\tilde{\phi}^s_{\omega}(\xi+\epsilon\eta)-\tilde{\phi}^s_{\omega}(\xi)) \, \d \theta \|^{2} \, \mathrm{d}s\right)^{\frac{1}{2}}\\
			&\leq Me^{\omega_A  t}\|\eta\|_{H_0}+C_{\beta} k(r) \left( \int_0^t\|P^s(\epsilon,\eta,\omega)\|^{2} \, \mathrm{d}s\right)^{\frac{1}{2}}\\
			&\leq Me^{\omega_A  t+C_{\beta}k(r)t}\|\eta\|_{H_0}.
		\end{align*}
		Using the dominated convergence theorem and letting $\epsilon\to 0$, we find that
		\begin{equation*}
			D_{\xi}\tilde{\phi}^t_{\omega}[\eta]=T_{0}(t)\eta+\int_0^t (-A_0)^{\beta}T_{0}(t-s)(-A)^{-\beta}DG_{\tilde{\phi}^s_{\omega}(\xi)}(D_{\xi}\tilde{\phi}^t_{\omega}[\eta]) \, \mathrm{d}s.
		\end{equation*}
		This completes the proof.
	\end{proof}

	In the result that follows, we achieve an integrable bound for our solutions to the random PDE \eqref{eqn:random_PDE}.
	\begin{theorem}\label{Priory-bound-cocycle}
		Let assumptions \ref{Assumption1.1} and \ref{1.2} be satisfied and assume that $\omega_A<0$. Assume that $G:H_0\to H$ is of class $\mathcal{C}^1$. Then the random dynamical system $\tilde{\phi}:\mathbb{R^+}\times\Omega\times H_0\to H_0$ defined on \eqref{cocycle-NP} satisfies the following apriori bound
		\begin{align*}
			\Vert\tilde{\phi}^{t}_{\omega}(\xi)\Vert\leq \kappa(t)\Vert\xi\Vert+b(t,\omega)
		\end{align*}
		where $\kappa:(0,\infty)\rightarrow (0,\infty)$ is a continuous function and $b:(0,\infty)\times\Omega\rightarrow (0,\infty)$ is a random variable satisfying
		\begin{align*}
			\sup_{0\leq t\leq 1}b(t,\cdot)\in \bigcap_{1 \leq p < \infty}{L}^{p}(\Omega).
		\end{align*}
	\end{theorem}
	\begin{proof}
		As pointed it out in Theorem \ref{Existence-solution}, the mild solution of the random PDE \eqref{eqn:random_PDE} satisfies, for any $\xi\in H_0$
		\begin{align*}
			V(t,\omega)=T_{0}(t)(\xi-Y_\omega)+\int_{0}^{t}(-A_0)^{\beta}T_{0}(t-s)(-A)^{-\beta} G\big(V(s,\omega)+Y_{\theta_{s}\omega})\big) \, \mathrm{d}s.
		\end{align*}
		As $G$ is Fréchet differentiable, it follows that $G$ has a linear growth, i.e., there exist constants $\kappa_1,\kappa_2>0$ such that for any $h\in H_0$,
		\begin{align*}
			\Vert G(h)\Vert\leq\kappa_1+\kappa_2\Vert h\Vert.
		\end{align*}
		Now, by combining the estimates \eqref{crucial-estimate1}, \eqref{Crucial-equlity-isg} and \eqref{Crucial-inequlity-isg}, we obtain for any $\beta\in(1-\frac{1}{p^{*}},1)$
		\begin{align*}
			\Vert V(t,\omega)\Vert&\leq \Vert T_0(t)(\xi-Y_\omega)\Vert+\int_{0}^{t}\big\Vert(-A_0)^{\beta}T_{0}(t-s)(-A)^{-\beta} G\big(V(s,\omega) + Y_{\theta_{s}\omega})\big)\big\Vert \, \mathrm{d}s\\
			&\leq \Vert T_{0}(t)(\xi-Y_\omega)\Vert\\
			&\quad + M_{\beta}\Vert(-A)^{-\beta}\Vert\int_{0}^{t}(t-s)^{-\beta}\exp(\omega_A(t-s))\big( \kappa_1+\kappa_2\Vert V(s,\omega)+Y_{\theta_{s}\omega}\Vert\big) \, \mathrm{d}s \\
			&\leq \underbrace{\Vert T_{0}(t)(\xi-Y_\omega)\Vert+M_{\beta}\Vert(-A)^{-\beta}\Vert\int_{0}^{t}(t-s)^{-\beta}\exp(\omega_A(t-s))\big( \kappa_1+\kappa_2\Vert Y_{\theta_{s}\omega}\Vert \big) \, \mathrm{d}s}_{ \eqqcolon a(t,\omega)} \\
			&\quad + \kappa_2M_{\beta}\Vert(-A)^{-\beta}\Vert\int_{0}^{t}(t-s)^{-\beta}\exp(\omega_A(t-s))\Vert V(s,\omega)\Vert \, \mathrm{d}s.
		\end{align*}
		Therefore, we obtain for any $\omega_A<0$ and $t>0$
		\begin{align*}
			\exp(-\omega_{A}t) \Vert V(t,\omega)\Vert \leq \exp(-\omega_{A}t)a(t,\omega)+\kappa_2M_{\beta}\Vert(-A)^{-\beta}\Vert\int_{0}^{t}(t-s)^{-\beta}\exp(-\omega_As)\Vert V(s,\omega)\Vert\mathrm{d}s.
		\end{align*}
		We infer from Lemma \ref{powered_Gronwall-type_inequality} and Lemma \ref{sum_bound} 
		\begin{align}\label{BVBVC}
			\begin{split}
			&\Vert V(t,\omega)\Vert \\ 
			\leq\ &a(t,\omega)+\int_{0}^{t}\sum_{n\geq 1}\frac{\big(\kappa_2M_{\beta}\Vert(-A)^{-\beta}\Gamma(1-\beta)\big)^n}{\Gamma(n(1-\beta))}(t-s)^{n(1-\beta)-1}\exp(\omega_{A}(t-s))a(s,\omega) \, \mathrm{d}s\\
			\leq &a(t,\omega) +\kappa_2M_{\beta}\Vert(-A)^{-\beta}\Vert\Gamma(1-\beta)\int_{0}^{t} R\big(\kappa_2M_{\beta}\Vert(-A)^{-\beta}\Vert\Gamma(1-\beta)(t-s)\big)(t-s)^{-n\beta}\exp(\omega_{A}(t-s))a(s,\omega)\mathrm{d}s
						\end{split}
		\end{align}
		where $R$ is defined in \eqref{sum_bound_1}. Now, our claim follows from \eqref{BVBVC} and the fact that $\tilde{\phi}^{t}_{\omega}(\xi) = V(t,\omega) + Y_{\theta_t \omega}$.
	\end{proof}
	
	In the next proposition, we formulate a bound for the derivative of the cocycle map.
	\begin{proposition}\label{Derivative_growth}
		Let assumptions \ref{Assumption1.1} and \ref{1.2} hold with $\omega_A<0$ and let $G:H_0\to H$ be of class $\mathcal{C}^1$. Assume that there exists an increasing polynomial $p$ such that
		\begin{equation*}
			\|D_{\xi}G\|\leq p(\|\xi\|),\quad\forall\xi\in H_0.
		\end{equation*}
		Then there exist constants $M_1,M_2\geq 1$ such that
		\begin{align*}
			\left\Vert D_{\xi}\tilde{\phi}^t_{\omega}\right\Vert_{\mathcal{L}(H_{0})}
			\leq\|T_{0}(t)\|+M_1 P(t,\xi,\omega)^{\frac{1}{(1-\beta)}}t\exp(M_2P(t,\xi,\omega)^{\frac{1}{(1-\beta)}}t)=:\Gamma(t,\xi,\omega),
		\end{align*}
		where
		\begin{align*}
			P(t,\xi,\omega)=\sup_{0\leq s\leq t}p(\|\tilde{\phi}^s_{\omega}(\xi)\|).
		\end{align*}
	\end{proposition}
	\begin{proof}
		From Proposition \ref{Fréchet}, we have for any $\xi,\eta\in H_0$
		\begin{equation*}
			D_{\xi}\tilde{\phi}^t_{\omega}[\eta]=T_{0}(t)\eta+\int_0^t (-A_0)^{\beta}T_{0}(t-s)(-A)^{-\beta}DG_{\tilde{\phi}^s_{\omega}(\xi)}(D_{\xi}\tilde{\phi}^s_{\omega}[\eta]) \, \mathrm{d}s.
		\end{equation*}
		Therefore, we find that 
		\begin{align*}
			\|D_{\xi}\tilde{\phi}^t_{\omega}[\eta]\|\leq \|T_{0}(t)\eta\|+\int_0^t \|(-A_0)^{\beta}T_{0}(t-s)(-A)^{-\beta}DG_{\tilde{\phi}^s_{\omega}(\xi)}(D_{\xi}\tilde{\phi}^s_{\omega}[\eta])\|\mathrm{d}s.
		\end{align*}
		According to Lemma \ref{powered_Gronwall-type_inequality}, we obtain for any $t>0$ and $\xi,\eta\in H_0$
		\begin{align*}
			\exp{(-\omega_{A}t)}\|D_{\xi}\tilde{\phi}^t_{\omega}[\eta]\|&\leq \exp{(-\omega_{A}t)} \|T_{0}(t)\eta\|\\
			&\quad +M_{\beta}\|(-A)^{-\beta}\|\int_0^t (t-s)^{-\beta}\exp{(-\omega_{A}s)}p(\|\tilde{\phi}^s_{\omega}(\xi)\|)\|(D_{\xi}\tilde{\phi}^s_{\omega}[\eta])\|\, \mathrm{d}s\\
			&\leq\exp{(-\omega_{A}t)} \|T_{0}(t)\eta\|\\
			&\quad +\int_{0}^{t}\sum_{n\geq 1}\frac{\big[{g}(t,\omega)\Gamma(1-\beta)\big]^n}{\Gamma(n(1-\beta))}(t-s)^{n(1-\beta)-1}\exp(\omega_{A}(-s))\|T_{0}(s)\eta\| \, \mathrm{d}s
		\end{align*}
		where $g(t,\xi,\omega):=M_{\beta}\|(-A)^{-\beta}\|\sup_{0\leq s\leq t}p(\|\tilde{\phi}^s_{\omega}(\xi)\|)$. Again from Lemma \ref{sum_bound}, 
		\begin{comment}
			\begin{align}\label{TR}
				\begin{split}
					&\|D_{\xi}\tilde{\phi}^t_{\omega}[\eta]\|\leq \|T_{0}(t)\eta\|\int_{0}^{t}\sum_{n\geq 1}\frac{\big[{g}(t,\omega)\Gamma(1-\beta)\big]^n}{\Gamma(n(1-\beta))}(t-s)^{(n-1)(1-\beta)-\beta}\exp(\omega_{A}(t-s))\|T_{0}(s)\eta\|\mathrm{d}s\\
					& =  \|T_{0}(t)\eta\|+\int_{0}^{t}\sum_{n\geq 1}\frac{\big[{g}(t,\omega)\Gamma(1-\beta)(t-s)^{1-\beta}\big]^{n-1}}{\Gamma(n(1-\beta))}{g}(t,\omega)\Gamma(1-\beta)(t-s)^{-\beta}\exp(\omega_{A}(t-s))\|T_{0}(s)\eta\|\mathrm{d}s\\&\leq\|T_{0}(t)\eta\|+{g}(t,\omega)\Gamma(1-\beta)\int_{0}^{t} R\big({g}(t,\omega)\Gamma(1-\beta)(t-s)^{1-\beta}\big)(t-s)^{-\beta}\exp(\omega_{A}(t-s))\|T_{0}(s)\eta\|\mathrm{d}s.
				\end{split}
			\end{align}
			Here $R$ is defined on \eqref{sum_bound_1}. So \eqref{TR} yields our claim.
			$	\frac{\beta+1}{1-\beta}({g}(t,\omega)\Gamma(1-\beta)(t-s)^{1-\beta})^{\frac{1+\beta}{1-\beta}}+\frac{\alpha^{\frac{2}{1-\beta}}\exp(\alpha^{\frac{1}{1-\beta}})}{1-\beta}$
		\end{comment}
		there is an $M_A\geq1$ such that
		\begin{align*}
			\|D_{\xi}\phi^t_{\omega}[\eta]\| &\leq \|T_{0}(t)\eta\|\\
			&\quad +\int_{0}^{t}\sum_{n\geq 1}\frac{\big[{g}(t,\xi,\omega)\Gamma(1-\beta)\big]^n}{\Gamma(n(1-\beta))}(t-s)^{n(1-\beta)-1}\exp(\omega_{A}(t-s))\|T_{0}(s)\eta\| \, \mathrm{d}s\\
			&\leq \|T_{0}(t)\eta\|\\
			&\quad + \sum_{n\geq 1}\frac{\big[{g}(t,\xi,\omega)\Gamma(1-\beta)\big]^n}{\Gamma(n(1-\beta))}\int_{0}^{t}(t-s)^{n(1-\beta)-1}\exp(\omega_{A}(t-s))\|T_{0}(s)\eta\| \, \mathrm{d}s\\
			&\leq \|T_{0}(t)\eta\| + M_A \sum_{n\geq 1}\frac{\big[{g}(t,\xi,\omega)\Gamma(1-\beta)t^{(1-\beta)}\big]^n}{n(1-\beta)\Gamma(n(1-\beta))}\|\eta\|_{H_0}\\
			&\leq \|T_{0}(t)\eta\| + M_A \sum_{n\geq 1}\frac{\big[{g}(t,\xi,\omega)^{\frac{1}{(1-\beta)}}\Gamma(1-\beta)^{\frac{1}{(1-\beta)}}t\big]^{n(1-\beta)}}{\lfloor n(1-\beta)! \rfloor}\|\eta\|_{H_0}\\
			&\leq \|T_{0}(t)\eta\| + M_A {g}(t,\xi,\omega)^{\frac{1}{(1-\beta)}}\Gamma(1-\beta)^{\frac{1}{(1-\beta)}}t \sum_{n\geq 1}\frac{\big[{g}(t,\omega)^{\frac{1}{(1-\beta)}}\Gamma(1-\beta)^{\frac{1}{(1-\beta)}}t\big]^{n(1-\beta)-1}}{\lfloor n(1-\beta)! \rfloor}\|\eta\|_{H_0}\\
			&\leq \|T_{0}(t)\eta\| + M_A {g}(t,\xi,\omega)^{\frac{1}{(1-\beta)}}\Gamma(1-\beta)^{\frac{1}{(1-\beta)}}t \sum_{n\geq 0}\frac{\big[{g}(t,\omega)^{\frac{1}{(1-\beta)}}\Gamma(1-\beta)^{\frac{1}{(1-\beta)}}t\big]^{\lfloor n(1-\beta)\rfloor}}{\lfloor n(1-\beta)! \rfloor}\|\eta\|_{H_0}\\
			&\leq\|T_{0}(t)\eta\| + M_A {g}(t,\xi,\omega)^{\frac{1}{(1-\beta)}}\Gamma(1-\beta)^{\frac{1}{(1-\beta)}}t\exp{({g}(t,\omega)^{\frac{1}{(1-\beta)}}\Gamma(1-\beta)^{\frac{1}{(1-\beta)}}t)}\|\eta\|_{H_0}.
		\end{align*}
		From Theorem \ref{Priory-bound-cocycle}, we know that
		\begin{equation*}
			g(t,\xi,\omega)\leq M_{\beta}\|(-A)^{-\beta}\|\sup_{0\leq s\leq t}p(\|\kappa(s)\Vert\xi\Vert+b(s,\omega)\|.
		\end{equation*}
		This completes the proof.
	\end{proof}
	
	\begin{corollary}\label{estimate on differnce}
		Let assumptions \ref{Assumption1.1} and \ref{1.2} be satisfied with $\omega_A<0$, and let $G:H_0\to H$ be of class $\mathcal{C}^1$. Assume that there exists an increasing polynomial $p$ such that
		\begin{equation*}
			\|D_{\xi}G\|\leq p(\|\xi\|),\quad\forall\xi\in H_0.
		\end{equation*}
		Then for any initial states $\xi_1,\xi_2\in H_0$
		\begin{align*}
			\Vert\tilde{\phi}_{\omega}^{t}(\xi_2)-\tilde{\phi}^{t}_{\omega}(\xi_1)\Vert
			&\leq\big(\|T_{0}(t)\|+M_1 \tilde{P}(t,\xi_2,\xi_1,\omega)^{\frac{1}{(1-\beta)}}t\exp(M_2\tilde{P}(t,\xi_2,\xi_1,\omega)^{\frac{1}{(1-\beta)}}t)\big)\|\xi_2-\xi_1\|\\
			&\quad =:\tilde\Gamma(t,\xi_2,\xi_1,\omega)\|\xi_2-\xi_1\|,
		\end{align*}
		where
		\begin{align*}
			\tilde{P}(t,\xi_2,\xi_1,\omega):=\sup_{0\leq s\leq t}p\big(\sup_{0\leq\theta\leq1}\|\tilde{\phi}^s_{\omega}(\theta\xi_2+(1-\theta)\xi_1)\|\big)
		\end{align*}
	\end{corollary}
	
	\begin{proof}
		We use the estimate provided in Proposition \ref{Derivative_growth} since
		\begin{align*}
			\Vert\tilde{\phi}_{\omega}^{t}(\xi_2)-\tilde{\phi}^{t}_{\omega}(\xi_1)\Vert\leq\int_{0}^{1}\Vert D_{\theta\xi_2+(1-\theta)\xi_1}\tilde{\phi}^{t}_{\omega}\Vert\Vert\xi_2-\xi_1\Vert\mathrm{d}\theta,
		\end{align*}
		and the statement follows directly.
	\end{proof}
We will now formulate an additional technical result that will be necessary for our subsequent analysis.
	\begin{proposition}\label{PPR}
		Let assumptions \ref{Assumption1.1} and \ref{Assumption1.2} be satisfied. Assume that $G$ is Fréchet differentiable and that there is some $0<r\leq 1$ such that for any $\xi_2,\xi_1\in H_0$,
		\begin{align*}
			\Vert D_{\xi_2}G-D_{\xi_1}G\Vert_{\mathcal{L}(H_0)}\leq \Vert\xi_2-\xi_1\Vert^rQ(\Vert\xi_2\Vert,\Vert\xi_1\Vert),
		\end{align*}
		where $Q$ is an increasing polynomial. Then for any  $\xi_1,\xi_2$ and $\eta$ in $H_0$,
		\begin{align*}
			\sup_{0\leq t\leq 1}\Vert D_{\xi_2}\phi^t_{\omega}[\eta]-D_{\xi_1}\phi^t_{\omega}[\eta]\Vert\leq  \exp\big(Q(\Vert\xi_1\Vert,\Vert\xi_2\Vert,f(\omega))\big)\Vert\xi_2-\xi_1\Vert^{r}\Vert\eta\Vert		
		\end{align*}
		where $f(\omega)\in\bigcap_{1 \leq p < \infty}L^{p}(\Omega).$
	\end{proposition}
	
	\begin{proof}
		Let $\xi_1,\xi_2,\eta\in H_0$ and $t>0$. According to Proposition \ref{Fréchet}, we have
		\begin{align*}
			&D_{\xi_2}\tilde{\phi}^t_{\omega}[\eta]-D_{\xi_1}\tilde{\phi}^t_{\omega}[\eta]\\&\quad=\int_{0}^{t}(-A_0)^{\beta}T_{0}(t-s)(-A)^{-\beta}\big(D_{\tilde{\phi}^s_{\omega}(\xi_2)}G(D_{\xi_2}\tilde{\phi}^s_{\omega}[\eta])-D_{\tilde{\phi}^s_{\omega}(\xi_1)}G(D_{\xi_1}\tilde{\phi}^s_{\omega}[\eta])\big) \, \mathrm{d}s.
		\end{align*}
		If we put $$L(s,\xi_2,\xi_1,\eta):=D_{\xi_2}\tilde{\phi}^s_{\omega}[\eta]-D_{\xi_1}\tilde{\phi}^s_{\omega}[\eta].$$ Then
		\begin{align*}
			&\Vert L(t,\xi_2,\xi_1,\eta)\Vert \\&\quad \leq  \int_{0}^{t}\bigg\Vert(-A_0)^{\beta}T_{0}(t-s)(-A)^{-\beta}\big(D_{\tilde{\phi}^s_{\omega}(\xi_2)}G-D_{\tilde{\phi}^s_{\omega}(\xi_1)}G\big)(D_{\xi_2}\tilde{\phi}^s_{\omega}[\eta])\bigg\Vert\mathrm{d}s\\&\qquad+\int_{0}^{t} \bigg\Vert(-A_0)^{\beta}T_{0}(t-s)(-A)^{-\beta}D_{\tilde{\phi}^s_{\omega}(\xi_1)}G\big(L(s,\xi_2,\xi_1,\eta)\big)\bigg\Vert \, \mathrm{d}s
			\\ &\quad\leq h(t)+M_{\beta}\Vert (-A)^{-\beta}\Vert\sup_{0\leq s\leq t}\Vert D_{\tilde{\phi}^s_{\omega}(\xi_1)}G\Vert_{\mathcal{L}(H_0)}\int_{0}^{t}(t-s)^{-\beta}\exp(\omega_{A}(t-s))\Vert L(s,\xi_2,\xi_1,\eta)\Vert \, \mathrm{d}s
		\end{align*}
		where
		\begin{align*}
			h(t):=\int_{0}^{t}\bigg\Vert(-A_0)^{\beta}T_{0}(t-s)(-A)^{-\beta}\big(D_{\tilde{\phi}^s_{\omega}(\xi_2)}G-D_{\tilde{\phi}^s_{\omega}(\xi_1)}G\big)(D_{\xi_2}\tilde{\phi}^s_{\omega}[\eta])\bigg\Vert \, \mathrm{d}s.
		\end{align*}
		
		It follows that
		\begin{align*}
			\exp(-\omega_At)\Vert L(t,\xi_2,\xi_1,\eta)\Vert\leq \exp(-\omega_At) h(t)+  \tilde{\kappa}(t)\int_{0}^{t}(t-s)^{-\beta}\exp(-\omega_{A}s)\Vert L(s,\xi_2,\xi_1,\eta)\Vert \, \mathrm{d}s,
		\end{align*}
		where
		\begin{align*}
			\tilde{\kappa}(t):=M_{\beta}\Vert (-A)^{-\beta}\Vert\sup_{0\leq s\leq t}\Vert D_{\tilde{\phi}^s_{\omega}(\xi_1)}G\Vert_{\mathcal{L}(H_0)}\exp(-\omega_At).
		\end{align*} 
		According to Lemma \ref{powered_Gronwall-type_inequality},
		\begin{align}\label{L}
			\begin{split}    
				& \Vert L(t,\xi_2,\xi_1,\eta)\Vert\leq \\&\quad h(t)+\int_{0}^{t}\sum_{n\geq 1}\frac{[\tilde{\kappa}(t)\Gamma(1-\beta)]^n}{\Gamma(n(1-\beta))}(t-s)^{n(1-\beta)-1}\exp(\omega_A(t-s))h(s) \, \mathrm{d}s.
			\end{split}
		\end{align}
		We proceed by estimating $h(t)$. For any $t>0$, we have
		\begin{align*}
			h(t) &\leq M_{\beta}\Vert (-A)^{-\beta}\Vert \int_{0}^{t}(t-s)^{-\beta}\exp(\omega_{A}(t-s))\big\Vert\big(D_{\tilde{\phi}^s_{\omega}(\xi_2)}G-D_{\tilde{\phi}^s_{\omega}(\xi_1)}G\big)(D_{\xi_2}\tilde{\phi}^s_{\omega}[\eta])\big\Vert \, \mathrm{d}s\\
			&\leq M_{\beta}\Vert (-A)^{-\beta}\Vert\int_{0}^{t}(t-s)^{-\beta}\exp(\omega_{A}(t-s))\Vert\tilde{\phi}^s_{\omega}(\xi_2)-\tilde{\phi}^s_{\omega}(\xi_1)\Vert^r Q\big(\Vert \phi^s_{\omega}(\xi_2)\Vert,\Vert \tilde{\phi}^s_{\omega}(\xi_1)\Vert\big)\Vert D_{\xi_2}\tilde{\phi}^s_{\omega}[\eta]\Vert \, \mathrm{d}s.
		\end{align*}
		From Proposition \ref{Derivative_growth} and Corollary \ref{estimate on differnce},
		\begin{align}\label{h_}
			\begin{split}    
				h(t) &\leq M_{\beta}\Vert (-A)^{-\beta}\int_{0}^{t}s^{-\beta}\exp(\omega_{A}s)\mathrm{d}s \Vert\sup_{0\leq s\leq t}\bigg[Q\big(\Vert\tilde{\phi}^s_{\omega}(\xi_2)\Vert,\Vert \tilde{\phi}^s_{\omega}(\xi_1)\Vert\big)\tilde\Gamma(s,\xi_2,\xi_1,\omega)^r\Gamma(s,\xi_2,\omega)\bigg]\\ 
				&\qquad \times\Vert\xi_2-\xi_1\Vert^{r}\Vert\eta\Vert \\
				&=: \Lambda(t,\xi_2,\xi_1,\omega) \Vert\xi_2-\xi_1\Vert^{r}\Vert\eta\Vert
			\end{split}
		\end{align}
		where $\Lambda$ is implicitly defined and $\Gamma$ and $\tilde\Gamma$ are taken as in Proposition \ref{Derivative_growth} and Corollary \ref{estimate on differnce}. We use inequality \eqref{h_} to estimate $\Vert L(t,\xi_2,\xi_1,\eta)\Vert$ in \eqref{L}. Since by definition, $\Lambda(t,\xi_2,\xi_1,\omega)$ is an increasing function in the first argument, we have
		\begin{align}\label{ASW}
			\begin{split}
				\Vert L(t,\xi_2,\xi_1,\eta)\Vert&\leq \Lambda(t,\xi_2,\xi_1,\omega) \Vert\xi_2-\xi_1\Vert^{r}\Vert\eta\Vert\\
				&\quad\times \left[ 1+ \int_{0}^{t}\sum_{n\geq 1}\frac{(\tilde{\kappa}(t)\Gamma(1-\beta))^n}{\Gamma(n(1-\beta))}(t-s)^{n(1-\beta)-1}\exp(\omega_A(t-s)) \, \mathrm{d}s \right].
			\end{split}
		\end{align}
		To finish the proof, we must find a bound for 
		\begin{align*}
			\Lambda(t,\xi_2,\xi_1,\omega) \left[ 1+ \int_{0}^{t}\sum_{n\geq 1}\frac{[\tilde{\kappa}(t)\Gamma(1-\beta)]^n}{\Gamma(n(1-\beta))}(t-s)^{n(1-\beta)-1}\exp(\omega_A(t-s)) \, \mathrm{d}s \right].
		\end{align*}
		From Lemma \eqref{powered_Gronwall-type_inequality},
		\begin{align}\label{Bound_3}
			\begin{split}
				&\int_{0}^{t}\sum_{n\geq 1}\frac{[\tilde{\kappa}(t)\Gamma(1-\beta)]^n}{\Gamma(n(1-\beta))}(t-s)^{n(1-\beta)-1}\exp(\omega_A(t-s)) \, \mathrm{d}s\\
				&\quad \leq\int_{0}^{t}\tilde{\kappa}(t)\Gamma(1-\beta)(t-s)^{-\beta}\exp(\omega_A(t-s))R\big(\tilde{\kappa}(t)(t-s)^{1-\beta}\Gamma(1-\beta)\big) \, \mathrm{d}s\\
				&\qquad \leq  \tilde{\kappa}(t)\Gamma(1-\beta)R\big(\tilde{\kappa}(t)t^{1-\beta}\Gamma(1-\beta)\big)\int_{0}^{t}(t-s)^{-\beta}\exp(\omega_A(t-s)) \, \mathrm{d}s.
			\end{split}
		\end{align}
		Recall the definition of $\tilde{\kappa}(t)$ and of the function $R$ in \eqref{sum_bound_1}. We can then conclude that
		\begin{align*}
			&\Vert L(t,\xi_2,\xi_1,\eta)\Vert \\
			\leq\ &\Vert\xi_2-\xi_1\Vert^{r}\Vert\eta\Vert \Lambda(t,\xi_2,\xi_1,\omega) \left[ \tilde{\kappa}(t)\Gamma(1-\beta)R\big(\tilde{\kappa}(t)t^{1-\beta}\Gamma(1-\beta)\big)\int_{0}^{t}(t-s)^{-\beta}\exp(\omega_A(t-s))\mathrm{d}s\right]
		\end{align*}
		To finalize the proof, it is enough to use \eqref{Bound_3} and the results of Theorem \ref{Priory-bound-cocycle}, Proposition \ref{Derivative_growth} and Corollary \ref{estimate on differnce}.
		%    \begin{align*}
			%  \exp(-\omega_At)\Vert L(t,\xi_2,\xi_1,\eta)\Vert\leq \exp(-\omega_At) h(t)+\int_{0}^{t}
			%  \end{align*}
	\end{proof}
	The following assumption enables us to prove the compactness of the linearized cocycle $D_{\xi}\phi^t_{\omega}:=\tilde{\psi}^{t_0}_{\omega}:H_{0}\rightarrow H_0$.
	\begin{assumption}\label{compactness-assum}
		Assume that $R(\lambda,A):=(\lambda-A)^{-1}$ is compact for $\lambda\in\rho(A)$.
	\end{assumption}
	\begin{remark}
		As established in \cite[Theorem 3.3]{Pazy-book}, the semigroup $T_0(t)$ is compact for $t>0$ if and only if it is continuous in the uniform operator topology for $t>0$ and if Assumption \ref{compactness-assum} holds. Given that $A_0$ is sectorial and Hille-Yosida, $T_0(t)$ is continuous in the uniform operator topology for $t>0$. Therefore, Assumption \ref{compactness-assum} is sufficient to ensure the compactness of the semigroup $T_0(t)$ for $t>0$.
	\end{remark}
	
	\begin{proposition}\label{compactness-result}
		Let assumptions \ref{Assumption1.1}, \ref{Assumption1.2} and \ref{compactness-assum} be satisfied. Assume that $\omega_A<0$ and $G:H_0\to H$ be of class $\mathcal{C}^1$. Then the linearized cocycle $\tilde{\psi}^{t}_{\omega} \coloneqq D_{\xi}\tilde{\phi}^t_{\omega}$ given by
		\begin{equation}\label{Fréchet_differentiable2-4}
			D_{\xi}\tilde{\phi}^t_{\omega}[\eta] = T_{0}(t)\eta + \lim_{\lambda\to\infty}\int_0^tT_{0}(t-s)\lambda(\lambda-A)^{-1}DG_{\tilde{\phi}^s_{\omega}(\xi)}(D_{\xi}\tilde{\phi}^s_{\omega}[\eta]) \, \mathrm{d}s, \quad \eta\in H_0,\quad t>0,
		\end{equation}
		is compact.
	\end{proposition}
	\begin{proof}
		To prove that $\tilde{\psi}^{t}_{\omega}$ is compact, we proceed by showing that the image of any bounded subset of $H_0$ under $\tilde{\psi}^{t}_{\omega}$ has a compact closure for $t>0$. We let $\mathcal{B}$ be a bounded subset of $H_0$, our goal is to demonstrate that the following set
		\begin{align*}
			\tilde{\psi}^{t}_{\omega}\mathcal{B} = \left\{\tilde{\psi}^{t}_{\omega}[\eta]~~ \,|\, ~~\eta \in \mathcal{B}\right\}
		\end{align*}
		has a compact closure for $t>0$.\\   
		Let $\epsilon>0$ and $\operatorname{Re}\lambda>\omega_A$. Then
		\begin{align*}
			\tilde{\psi}^{t}_{\omega}[\eta]=T_{0}(t)\eta+T_{0}(\epsilon)\lim_{\lambda\to\infty}\int_0^{t-\epsilon} T_{0}(t-s-\epsilon)&\lambda(\lambda-A)^{-1}DG_{\tilde{\phi}^s_{\omega}(\xi)}(D_{\xi}\tilde{\phi}^s_{\omega}[\eta])\mathrm{d}s\\+&\lim_{\lambda\to\infty}\int_{t-\epsilon}^t T_{0}(t-s)\lambda(\lambda-A)^{-1}DG_{\tilde{\phi}^s_{\omega}(\xi)}(D_{\xi}\tilde{\phi}^s_{\omega}[\eta])\mathrm{d}s.
		\end{align*}
		Clearly, if the assumption \ref{compactness-assum} holds, then the closure $\overline{T_0(t)\mathcal{B}_1}$ is compact for all bounded subsets $\mathcal{B}_1$ in $H_0$. Thus, if the following set 
		\begin{equation}\label{sufficient-conditi}
			\left\{\tilde{\psi}^{s}_{\omega}[\eta] \, \vert \,  \eta\in\mathcal{B},\quad s\in [0,t]\right\}.
		\end{equation}
		is bounded, then there exists a bounded subset $\mathcal{B}_1\subseteq H_0$ such that 
		\begin{equation*}
			K_\epsilon := \left\{\lim_{\lambda\to\infty}\int_0^{t-\epsilon} T_{0}(t-s-\epsilon)\lambda(\lambda-A)^{-1}DG_{\tilde{\phi}^s_{\omega}(\xi)}(D_{\xi}\tilde{\phi}^s_{\omega}[\eta]) \, \mathrm{d}s \, \vert \, \eta\in\mathcal{B}\right\}\subseteq\mathcal{B}_1.
		\end{equation*}
		This implies that $K_\epsilon$ has a compact closure in $H_0$. To verify this, we need to check whether the set \ref{sufficient-conditi} is bounded. Let $\eta\in\mathcal{B}$ and $s\in [0,t]$. According to Lemma \ref{powered_Gronwall-type_inequality} and using the estimates \eqref{Crucial-inequlity-isg} and \eqref{Derivative-of-G}, we find that
		\begin{align*}
			\|\tilde{\psi}^{s}_{\omega}[\eta]\|&\leq Me^{\omega_A s}\|\eta\|+M_{\beta}(\|(-A)^{-\beta}\|)
			\int_{0}^{s}(s-\sigma)^{-\beta}\exp(\omega_A(s-\sigma))DG_{\tilde{\phi}^\sigma_{\omega}(\xi)}(D_{\xi}\tilde{\phi}^\sigma_{\omega}[\eta]) \, \mathrm{d}\sigma\\
			&\leq M\int_{0}^{s}\sum_{n\geq 1}\frac{(\tilde{\kappa}(\beta)\Gamma(1-\beta))^n}{\Gamma(n(1-\beta))}(s-\sigma)^{n(1-\beta)-1}\exp(\omega_A(s-\sigma)) \, \mathrm{d}\sigma\|\eta\|,
		\end{align*}
		where $\tilde{\kappa}(\beta)=M_{\beta}(\|(-A)^{-\beta}\|)k(r)$. Therefore the set defined in \eqref{sufficient-conditi} is bounded. Consequently, $K_\epsilon$ is compact.\\
		We define the set
		\begin{equation*}
			\tilde{K}_\epsilon:=\left\{\lim_{\lambda\to\infty}\int_{t-\epsilon}^{t}T_{0}(t-s)\lambda(\lambda-A)^{-1}DG_{\tilde{\phi}^s_{\omega}(\xi)}(D_{\xi}\tilde{\phi}^s_{\omega}[\eta]) \, \mathrm{d}s \, \mid \,  \eta\in\mathcal{B}\right\}.
		\end{equation*}
		Using the same argument as before, every $y\in \tilde{K}_\epsilon$ satisfies $\|y\|\leq \epsilon C(\beta,M)\|\eta\|$ where $C(\beta,M)$ is a positive constant. If we choose $\epsilon > 0$ small enough that fits to domain of our definition, $\tilde{K}_\epsilon$ can be covered by one ball of radius $\frac{{\epsilon}}{2}$ centred at 0. So $\tilde{\psi}^t_{\omega}$ is relatively compact for $t>0$.
	\end{proof}
	\begin{remark}
Alternatively, compactness can be established by considering  
\[
(-{A}_0)^{\beta}{T}_{0}(t)(-A)^{-\beta} = (-{A}_0)^{-\epsilon}(-{A}_0)^{\beta+\epsilon}{T}_{0}(t)(-A)^{-\beta}
\]
for a sufficiently small $\epsilon > 0$ such that $\beta + \epsilon \in \left(1-\frac{1}{p^*}, 1\right)$. Since the operator  
\[
\int_{0}^{t}(-{A}_0)^{\beta + \epsilon}{T}_{0}(s)(-A)^{-\beta} \, \mathrm{d}s
\]
is bounded, and
\[
(-A_{0})^{-\epsilon} = \frac{1}{\Gamma(\epsilon)} \int_{0}^{\infty} s^{\epsilon-1} T_{0}(s) \, \mathrm{d}s,
\]
we can establish the compactness of $(-A_{0})^{-\epsilon}$. From this, using \eqref{Fréchet_differentiable}, we deduce the compactness of $\tilde{\psi}$.
	\end{remark}
	\section{Invariant Manifolds}
	In this section, we present our main results regarding the existence of invariant manifolds (stable, unstable, and center) around stationary points. These points can be understood as the random analogues to fixed points of a deterministic dynamical system. We will now give the formal definition.
	
	\begin{definition}\label{statinary point}
		Let $Z:\Omega\rightarrow H_{0}$ be a measurable function and let $\tilde{\phi}:[0,\infty)\times \Omega \times H_{0}\rightarrow H_{0}$ be the cocycle induced by the solutions to \eqref{evolution-pb}. We call $Z$ a \textbf{stationary point} for  $\tilde{\phi}$ if and only if
		\begin{align*}
			\tilde{\phi}^{t}_{\omega}(Z_{\omega})=Z_{\theta_{t}\omega},\quad \forall t\geq 0.		
		\end{align*}
		%		where $\Vert Z_{\omega}\Vert\in{\cap}_{p\geq 1}{L}^{p}(\Omega)$.
	\end{definition}
	The important property of a stationary point is that the linearized cocycle around it satisfies again the cocycle property, i.e. linearizing a non-linear cocycle around a stationary point yields a linear cocycle.\smallskip
	
	We will now formulate a version of the multiplicative ergodic theorem that holds for our case.

		\begin{theorem}[{Multiplicative ergodic theorem}]\label{MET}Let Assumption \ref{compactness-assum} hold, and assume that $G$ is Fréchet differentiable. Let $Z:\Omega\rightarrow H_0$ be a stationary point for the stochastic evolution equation \eqref{evolution-pb}. Set
			\begin{align}\label{eqn:lin_cocycle_around_Z}
				\tilde{\psi}^{t}_{\omega}(\cdot)\coloneqq D_{Z_{\omega}} \tilde{\phi}^t(\cdot).
			\end{align}
			Assume that for some $t_0>0$
			\begin{align}\label{V1}
				\max\big\lbrace\sup_{t\in [0,t_0]}\log^{+}\left\Vert \tilde{\psi}^{t}_{\omega}\right\Vert_{\mathcal{L}(H_{0})},\sup_{t\in [0,t_0]}\log^{+}\left\Vert\tilde{\psi}^{t_0-t}_{\theta_{t}\omega}\right\Vert_{\mathcal{L}(H_{0})}\big\rbrace\in{L}^{1}(\Omega).
			\end{align}
We define the set 
\[
F_{\mu}(\omega) := \left\{ \xi \in H_0 : \limsup_{t \rightarrow \infty} \frac{1}{t} \log \|\tilde{\psi}^{t}_{\omega}(\xi)\|_{\alpha} \leq \mu \right\},
\]
where \(\mu \in [-\infty, +\infty]\). Then, we have the Lyapunov exponents
\[
\mu_{1} > \mu_{2} > \ldots \in [-\infty, +\infty),
\]
which can either be finite and equal to \(-\infty\) for \(j \geq i\) or infinite and satisfying \(\lim_{n \rightarrow \infty} \mu_{n} = -\infty\). For every index \(i\) such that \(\mu_{i} > -\infty\), there exists a finite-dimensional subspace \(H^{i}_{\omega} \subset H_0\) satisfying the following properties:
			\begin{itemize}
				\item[(i)] (Invariance.)\ \ $\tilde{\psi}^{t}_{\omega}(H^{i}_{\omega})= H^i_{\theta_t \omega}$ for every $t \geq 0$.
				\item[(ii)] (Splitting.)\ \ $H_{\omega}^i \oplus F_{\mu_{i+1}}(\omega) = F_{\mu_i}(\omega)$ and $H_0 = \oplus_{1\leq j\leq i}H^j_{\omega}\oplus   F_{\mu_{i+1}}(\omega)$.
				\item[(iii)] (Fast growing subspace.)\  For each $ h\in H^{j}_{\omega} $,
				\begin{align*}
					&\lim_{t\rightarrow\infty}\frac{1}{t}\log\Vert \tilde{\psi}^{t}_{\omega}(h)\Vert = \mu_{j},\\
					& 		\lim_{t\rightarrow\infty}\frac{1}{t}\log\Vert (\tilde{\psi}^{t}_{\theta_{-t}\omega})^{-1}(h)\Vert =-\mu_{j}.
				\end{align*}
			\end{itemize}
		\end{theorem}
		\begin{proof}
		%	Let $t_0>0$. According to Proposition \ref{Derivative_growth} and Theorem \ref{Priory-bound-cocycle}, we have 
			%\begin{align}\label{V2}
			%	\sup_{t\in [0,t_0]}\log^{+}\left\Vert\tilde{\psi}^{t_0-t}_{\theta_{t}\omega}\right\Vert_{\mathcal{L}(H_{0})}\leq \sup_{t\in [0,t_0]}\Gamma(t-t_0,Y_{\theta_{t}\omega},\theta_t\omega)\in {L}^{1}(\Omega).
	%		\end{align}
			First, note that from the definitions of $Z_\omega$ and $\tilde{\psi}$, we can easily see that $\tilde{\psi}$ defines a linear cocycle. Recall that $H_0$ is a separable space, and thanks to our integrability condition \ref{V1}, we can apply \cite[Theorem 1.21]{GVR23} to obtain the stated results for the discrete-time cocycle $(\tilde{\psi}^{nt_0}_{\omega})_{(n,\omega) \in \mathbb{N} \times \Omega}$. We can generalize these results to continuous time by using the discrete-time version together with \eqref{V1}. The strategy for this is fairly standard, so we refer to \cite[Theorem 3.3 and Lemma 3.4]{LL10}, where it is explained in detail.
		\end{proof}
		
		\begin{notation} Assuming that the conditions of the MET hold, we introduce the following notation: Set $i_{0} \coloneqq \min \lbrace i \in \N \, \mid \,\mu_i<0 \rbrace$. If $\mu_1>0$, we define $j_{0} \coloneqq \max \lbrace i \in \N \, \mid\,  \mu_{i}>0 \rbrace$. Furthermore, we set 
			\begin{itemize}
				\item $S_\omega:=F_{\mu_{i_0}}(\omega)$.
				\item $U_{\omega}:=\oplus_{i: \mu_{i}>0}H^{i}_{\omega}$.
				\item $C_\omega:=H^{i_c}_{\omega}$ if there exists $\mu_{i_c}=0$.	
				\item $\Pi_{S_{\omega}}: S_{\omega}\oplus_{i:\mu_i=0}H^{i}_{\omega}\oplus U_{\omega}\rightarrow S_{\omega} $, projection into the $S_\omega$.
				\item  $\Pi_{U_{\omega}}: S_{\omega}\oplus_{i:\mu_i=0}H^{i}_{\omega}\oplus U_{\omega}\rightarrow U_\omega$, projection into the $U_\omega$.
			\end{itemize}
		\end{notation}

		The crucial condition to check in the MET is the integrability condition \eqref{V1}. In the next lemma, we formulate a sufficient condition in terms of the stationary point that implies this property.
		
		\begin{lemma}\label{integrability_derivative}
			Let assumptions \eqref{Assumption1.1} and \eqref{Assumption1.2} be satisfied and let $Z:\Omega\rightarrow H_0$ be a stationary point for the stochastic evolution equation \eqref{evolution-pb} that satisfies
			\begin{align}\label{eqn:moment_int_stat}
				\Vert Z_{\omega}\Vert\in \bigcap_{1 \leq p < \infty} {L}^{p}(\Omega).
			\end{align} 
			Assume that $G$ is Fréchet differentiable and there exists an increasing polynomial $p$ such that
			\begin{align}\label{condition on G}
				\Vert D_{Z_{\omega}}G\Vert_{\mathcal{L}(H_0,H)}\leq p(\Vert Z_\omega\Vert),\quad\forall \omega\in\Omega.  
			\end{align}
			Let $\tilde{\psi}^{t}_{\omega}$ be defined as in \eqref{eqn:lin_cocycle_around_Z}. Then it holds that
			\begin{align*}
		\max\big\lbrace\sup_{t\in [0,1]}\log^{+}\left\Vert \tilde{\psi}^{t}_{\omega}\right\Vert_{\mathcal{L}(H_{0})},\sup_{t\in [0,1]}\log^{+}\left\Vert\tilde{\psi}^{1-t}_{\theta_{t}\omega}\right\Vert_{\mathcal{L}(H_{0})}\big\rbrace\in{L}^{1}(\Omega).
			\end{align*}
		\end{lemma}
		
		\begin{proof}
	Let $t \geq 0$. According to Proposition \ref{Derivative_growth}, we have

\begin{align*}
    \left\Vert \tilde{\psi}^{t}_{\omega} \right\Vert_{\mathcal{L}(H_{0})} \leq \|T_{0}(t)\| + M_1 P(t, Z_\omega, \omega)^{\frac{1}{1 - \beta}} t \exp\left(M_2 P(t, Z_\omega, \omega)^{\frac{1}{1 - \beta}} t\right),
\end{align*}
where
\begin{align*}
    P(t, Z_\omega, \omega) = \sup_{0 \leq s \leq t} p(\|\tilde{\phi}^s_{\omega}(Z_\omega)\|).
\end{align*}
Since $p$ is an increasing polynomial, to check that $\sup_{t \in [0, 1]} \log^{+}\left\Vert \tilde{\psi}^{t}_{\omega} \right\Vert_{\mathcal{L}(H_{0})} \in L^{1}(\Omega)$, it is sufficient to verify that
\begin{align}\label{Integrability_2}
    \sup_{0 \leq t \leq 1} \|\tilde{\phi}^t_{\omega}(Z_\omega)\| \in \bigcap_{1 \leq p < \infty} L^{p}(\Omega).
\end{align}
This follows immediately from Theorem \ref{Priory-bound-cocycle}. Similarly, we can argue to obtain
\begin{align*}
   \sup_{0 \leq t \leq 1}  \left\Vert\tilde{\psi}^{1-t}_{\theta_{t}\omega}\right\Vert_{\mathcal{L}(H_{0})} \in L^{1}(\Omega).
\end{align*} 
		\end{proof}
		We can now formulate our main result about the existence of local stable manifolds. It is basically a reformulation of the abstract stable manifold theorem proven in \cite[Theorem 2.10]{GVR23}.  
		
		\begin{theorem}[Stable manifold theorem]\label{stable_manifold}
			Let assumptions \eqref{Assumption1.1}, \eqref{Assumption1.2} and \ref{compactness-assum} hold and let $Z:\Omega\rightarrow H_0$ be a stationary point for the stochastic evolution equation \eqref{evolution-pb} that satisfies \eqref{eqn:moment_int_stat}. Assume that $G$ is Fréchet differentiable and that there is some $0<r\leq 1$ such that for any $\xi_2,\xi_1\in H_0$,
			\begin{align*}
				\Vert D_{\xi_2}G-D_{\xi_1}G\Vert_{\mathcal{L}(H_0)}\leq \Vert\xi_2-\xi_1\Vert^rQ(\Vert\xi_2\Vert,\Vert\xi_1\Vert),
			\end{align*}
			where $Q$ is an increasing polynomial. Let $(\mu_j)_{j \in \N}$ be the Lyapunov exponents provided by Theorem \ref{MET} and set $\mu_{j_0} \coloneqq \max \lbrace \mu_{j} \, \mid \, \mu_j<0 \rbrace$. 
			We fix an arbitrary time step $t_{0}>0$. For $ 0 < \upsilon < -\mu_{j_0} $, we can find a family of immersed submanifolds $S^{\upsilon}_{loc}(\omega)$ of $H_0$ and a set of full measure $\tilde{\Omega}$ of $\Omega$ such that
			\begin{enumerate}[(i)]
				\item There are random variables $ R_{1}^{\upsilon}(\omega), R_{2}^{\upsilon}(\omega)$, which are positive and finite on $\tilde{\Omega}$, and
				\begin{align}\label{eqn:rho_temp}
					\liminf_{p \to \infty} \frac{1}{p} \log R_i^{\upsilon}(\theta_{pt_0} \omega) \geq 0, \quad i = 1,2
				\end{align}
				such that
				\begin{align}\label{invari}
					\begin{split}
						\big{\lbrace} \xi \in H_0\, :\, \sup_{n\geq 0}\exp(nt_0\upsilon)\Vert\tilde{\phi}^{nt_0}_{\omega}(\xi)-Z_{\theta_{nt_0}\omega}\Vert &<R_{1}^{\upsilon}(\omega)\big{\rbrace}\subseteq S^{\upsilon}_{loc}(\omega)\\&\subseteq	\big{\lbrace}  \xi \in H_0\, :\, \sup_{n\geq 0}\exp(nt_0\upsilon)\Vert\tilde{\phi}^{nt_0}_{\omega}(\xi) - Z_{\theta_{nt_0}\omega}\Vert<R_{2}^{\upsilon}({\omega})\big{\rbrace}.
					\end{split}
				\end{align}
				
				\item For $ S^{\upsilon}_{loc}(\omega)$, 
				\begin{align*}
					T_{Z_{\omega}}S^{\upsilon}_{loc}(\omega) = S_{\omega}.
				\end{align*}
				
				\item There is nonnegative random variable $N(\omega)$ such that for $ n\geq N(\omega)$,
				\begin{align*}
					\tilde{\phi}^{nt_0}_{\omega}(S^{\upsilon}_{loc}(\omega))\subseteq S^{\upsilon}_{loc}(\theta_{nt_0}\omega).
				\end{align*}
				
				\item For $ 0<\upsilon_{1}\leq\upsilon_{2}< - \mu_{j_{0}} $,
				\begin{align*}
					S^{\upsilon_{2}}_{loc}(\omega)\subseteq S^{\upsilon_{1}}_{loc}(\omega).
				\end{align*}
				Also, for $n\geq N(\omega)$,
				\begin{align*}
					\tilde{\phi}^{nt_0}_{\omega}(S^{\upsilon_{1}}_{loc}(\omega))\subseteq S^{\upsilon_{2}}_{loc}(\theta_{nt_0} \omega)
				\end{align*}
				and consequently for $ \xi\in S^{\upsilon}_{loc}(\omega)$,
				\begin{align}\label{eqn:contr_char}
					\limsup_{n\rightarrow\infty}\frac{1}{n}\log\Vert\tilde{\phi}^{nt_0}_{\omega}(\xi) - Z_{\theta_{nt_0}\omega}\Vert\leq  t_0\mu_{j_{0}}.
				\end{align}
				
				\item 
				\begin{align*}
					\limsup_{n\rightarrow\infty}\frac{1}{n}\log\bigg{[}\sup\bigg{\lbrace}\frac{\Vert\tilde{\phi}^{nt_0}_{\omega}(\xi)-\tilde{\phi}^{nt_0}_{\omega}(\tilde{\xi})\Vert }{\Vert \xi-\tilde{\xi}\Vert},\ \ \xi \neq\tilde{\xi},\  \xi,\tilde{\xi}\in S^{\upsilon}_{loc}(\omega) \bigg{\rbrace}\bigg{]}\leq t_0\mu_{j_{0}}.
				\end{align*}
			\end{enumerate}
		\end{theorem}
		
		\begin{proof}
			As before, we set $\tilde{\psi}^{t}_{\omega}(\cdot):=D_{Z_{\omega}}\tilde{\phi}^{t}_{\omega}(\cdot)$. From Lemma \ref{integrability_derivative}, we have
			\begin{align*}
				\sup_{t\in [0,1]}\log^{+}\left\Vert \tilde{\psi}^{t}_{\omega}\right\Vert_{\mathcal{L}(H_{0})}\in{L}^{1}(\Omega),
			\end{align*}
			thus the conditions of the MET \ref{MET} are satisfied. The idea is now to apply the abstract result stated in \cite[Theorem 2.10]{GVR23}. To establish this, first observe that the measurability condition \cite[Assumption 1.1]{GVR23} is clearly satisfied. This follows from the fact that $H_0$ is a closed subspace of the separable Hilbert space $H$, and is thus also separable. Moreover, for every $x, y \in H_0$, $t_0 > 0$, and $n, k \in \mathbb{N}$, the map
\[
\omega \mapsto \left\Vert \tilde{\psi}_{\theta_{n t_0} \omega}^{k t_0} \left( \tilde{\psi}_{\omega}^{n t_0} - y \right) \right\Vert
\]
is measurable. It remains to check that a bound of the form \cite[Equation (2.5)]{GVR23} holds in our case. We define
			\begin{align*}
				K_{\omega}(\xi) \coloneqq \tilde{\phi}^{t_0}_{\omega}(Z_{\omega}+\xi)-\tilde{\phi}^{t_0}_{\omega}(Z_{\omega})-\tilde{\psi}^{t_0}_{\omega}(\xi). 
			\end{align*}
		From Proposition \ref{PPR}, we obtain
		\begin{align}\label{IOPL}
			\begin{split}
				\Vert K_{\omega}(\xi) - K_{\omega}(\tilde\xi) \Vert
				&\leq \int_{0}^{1} \Vert \big( D_{Z_{\omega} + \theta\xi + (1-\theta)\tilde\xi} \tilde{\phi}^{t_0}_{\omega} - D_{Z_{\omega}} \tilde{\phi}^{t_0}_{\omega} \big) [\xi - \tilde{\xi}] \Vert \, \mathrm{d}\theta \\
				&\leq \exp\big(Q(\Vert \xi \Vert + \Vert \tilde\xi \Vert + \Vert Z_\omega \Vert, \Vert Z_\omega \Vert, f(\omega))\big) \Vert \xi - \tilde\xi \Vert (\Vert \xi \Vert^r + \Vert \tilde\xi \Vert^r),
			\end{split}
		\end{align}
		for any $\xi, \tilde{\xi} \in H_0$ and $\omega \in \Omega$, where $Q$ is an increasing polynomial. Note that we can find two increasing polynomials, $Q_1$ and $Q_2$, such that
		\begin{align*}
			Q(\Vert \xi \Vert + \Vert \tilde\xi \Vert + \Vert Z_\omega \Vert, \Vert Z_\omega \Vert, f(\omega)) \leq Q_1(\Vert \xi \Vert + \Vert \tilde\xi \Vert) + Q_2(\Vert Z_\omega \Vert, f(\omega))
		\end{align*}
		for any $\xi, \tilde{\xi} \in H_0$ and $\omega \in \Omega$. Therefore,
		\begin{align*}
			\Vert K_{\omega}(\xi) - K_{\omega}(\tilde\xi)\Vert
			&\leq \exp\big(Q_1(\Vert \xi \Vert + \Vert \tilde\xi \Vert)\big)\underbrace{\exp\big(Q_2(\Vert Z_\omega \Vert, f(\omega))\big)}_{R(\omega)} \Vert \xi - \tilde\xi \Vert (\Vert \xi \Vert^r + \Vert \tilde\xi \Vert^r) .
		\end{align*}
		Recall that $\Vert Z_\omega \Vert, f(\omega) \in \bigcap_{1 \leq p < \infty} L^p(\Omega)$. Since $Q_2$ is a polynomial, we have
		\begin{align*}
			\log^+ (R(\omega)) \in L^1(\Omega).
		\end{align*}
		Thus, from Birkhoff's Ergodic Theorem, on a set of full measure:
		\begin{align}
			\lim_{n \to \infty} \frac{1}{n} \log^+ (R(\theta_{nt_0} \omega)) = 0.
		\end{align}
		So, we have verified the necessary condition in \cite[Theorem 2.10]{GVR23}, and we can apply this result to prove our claims.
		\end{proof}

		Next, we formulate the unstable manifold theorem.
		
		\begin{theorem}[Unstable manifold theorem]\label{unstable_manifold}
			We impose the same assumptions as in Theorem \ref{stable_manifold}. Furthermore, we assume that $\mu_1 > 0$ and set $\mu_{j_0} \coloneqq \min \lbrace \mu_{j}\, \mid \, \mu_j > 0 \rbrace$. Fix an arbitrary time step $t_{0}>0$. For any $t \in \R$, set $\sigma_t \coloneqq \theta_{-t}$. Then, for any $ 0 < \upsilon < \mu_{j_0} $, there exists a family of immersed submanifolds $U^{\upsilon}_{loc}(\omega)$ of $H_0$ and a set of full measure $\tilde{\Omega}$ of $\Omega$ such that
			\begin{enumerate}[(i)]
				\item There are random variables $ \tilde{R}_{1}^{\upsilon}(\omega), \tilde{R}_{2}^{\upsilon}(\omega)$, which are positive and finite on $\tilde{\Omega}$, and
				\begin{align}\label{eqn:rho_temp2}
					\liminf_{p \to \infty} \frac{1}{p} \log\tilde{R}_i^{\upsilon}(\sigma_{pt_0} \omega) \geq 0, \quad i = 1,2
				\end{align}
				such that
				\begin{align}\label{invari2}
					\begin{split}
						&\bigg{\lbrace} \xi_{\omega} \in H_0 \, \mid \, \exists \lbrace \xi_{\sigma_{nt_{0}} \omega} \rbrace_{n\geq 1} \text{ s.t. }
						\tilde{\phi}^{mt_{0}}_{\sigma_{nt_{0}}\omega}(\xi_{\sigma_{nt_{0}}\omega}) = \xi_{\sigma_{(n-m)t_{0}}\omega} \text{ for all } 0 \leq m \leq n \text{ and } \\ 
						&\quad \sup_{n\geq 0} \exp(nt_0\upsilon) \| \xi_{\sigma_{nt_{0}} \omega} - Z_{\sigma_{nt_{0}}\omega} \| < \tilde{R}_{1}^{\upsilon}(\omega)\bigg{\rbrace} \subseteq U^{\upsilon}_{loc}(\omega)
						\subseteq \bigg{\lbrace} \xi_{\omega} \in H_0 \, \mid \, \exists \lbrace \xi_{\sigma_{nt_{0}} \omega} \rbrace_{n\geq 1} \text{ s.t. }\\
						&\tilde{\phi}^{mt_{0}}_{\sigma_{nt_{0}}\omega}(\xi_{\sigma_{nt_{0}}\omega}) = \xi_{\sigma_{(n-m)t_{0}}\omega} \text{ for all } 0 \leq m \leq n \text{ and } \sup_{n\geq 0} \exp(nt_0\upsilon) \| \xi_{\sigma_{nt_{0}} \omega} - Z_{\sigma_{nt_{0}}\omega} \| < \tilde{R}_{2}^{\upsilon}(\omega)\bigg{\rbrace}
					\end{split}
				\end{align}
				
				% \begin{align}\label{invari2}
					% 	\begin{split}
						% 		&\big{\lbrace} \xi_\omega \in H_0\, :\, \tilde{\phi}^{m t_0}_{\sigma_{nt_0}\omega}(\xi_{\sigma_{nt_0}\omega})=\xi_{\sigma_{(n-m)t_0}\omega}\,\text{ for all }\, 0\leq m\leq n \text{ and }\sup_{n\geq 0}\exp(nt_0\upsilon)\Vert \xi_{\sigma_{nt_0}\omega} - Z_{\sigma_{nt_0}\omega}\Vert <\tilde{R}_{1}^{\upsilon}(\omega)\big{\rbrace} \\ 
						%               \subseteq\ &U^{\upsilon}_{loc}(\omega) \\
						%               \subseteq\ &\big{\lbrace} \xi_\omega \in H_0\, :\, \tilde{\phi}^{m t_0}_{\sigma_{nt_0}\omega}(\xi_{\sigma_{nt_0}\omega})=\xi_{\sigma_{(n-m)t_0}\omega}\,\text{ for all }\, 0\leq m\leq n \text{ and }\sup_{n\geq 0}\exp(nt_0\upsilon)\Vert \xi_{\sigma_{nt_0}\omega} - Z_{\sigma_{nt_0}\omega}\Vert <\tilde{R}_{2}^{\upsilon}(\omega)\big{\rbrace}.
						% 	\end{split}
					%	\end{align}
				
				\item The submanifold $ U^{\upsilon}_{loc}(\omega)$ satisfies
				\begin{align*}
					T_{Z_{\omega}}U^{\upsilon}_{loc}(\omega) = U_\omega.
				\end{align*}
				
				\item For $ n\geq N(\omega) $,
				\begin{align*}
					U^{\upsilon}_{loc}(\omega)\subseteq\tilde{\phi}^{nt_0}_{\sigma_{nt_0}\omega}(U^{\upsilon}_{loc}(\sigma_{nt_0}\omega)) .
				\end{align*}
				
				\item For $ 0<\upsilon_{1}\leq\upsilon_{2}< \mu_{j_{0}} $,
				\begin{align*}
					U^{\upsilon_{2}}_{loc}(\omega)\subseteq U^{\upsilon_{1}}_{loc}(\omega).
				\end{align*}
				Also for $n\geq N(\omega) $,
				\begin{align*}
					U^{\upsilon_1}_{loc}(\omega)\subseteq\tilde{\phi}^{nt_0}_{\sigma_{nt_0}\omega}(U^{\upsilon_2}_{loc}(\sigma_{nt_0}\omega))
				\end{align*}
				and consequently, for $\xi_\omega \in U^{\upsilon}_{loc}(\omega)$,
				\begin{align}\label{eqn:contr_char2}
					\limsup_{n\rightarrow\infty}\frac{1}{n}\log\Vert \xi_{\sigma_{nt_0}\omega} - Z_{\sigma_{nt_0}\omega}\Vert\leq  -t_0\lambda_{j_{0}}.
				\end{align}
				
				\item 
				\begin{align*}
					\limsup_{n\rightarrow\infty}\frac{1}{n}\log\bigg{[}\sup\bigg{\lbrace}\frac{\Vert \xi_{\sigma_{nt_0}\omega} - \tilde{\xi}_{\sigma_{nt_0}\omega}\Vert }{\Vert \xi_{\omega} - \tilde{\xi}_\omega\Vert},\ \ \xi_{\omega} \neq\tilde{\xi}_\omega,\  \xi_{\omega},\tilde{\xi}_\omega\in U^{\upsilon}_{loc}(\omega) \bigg{\rbrace}\bigg{]}\leq -t_0\lambda_{\mu_{0}}.
				\end{align*}
			\end{enumerate}
		\end{theorem}
		\begin{proof}
			Using the same arguments as in the proof of Theorem \ref{stable_manifold}, the assertion follows from \cite[Theorem 2.17]{GVR23}.
		\end{proof}
		
		Finally, we state our result about the existence of center manifolds. Here, we will assume that one Lyapunov exponent is zero, i.e. $\mu_{i_c} = 0$. We set $\mu^{-} \coloneqq \mu_{i_c + 1}$. If $\mu_1 > 0$, we define $\mu^{+} \coloneqq \mu_{i_c - 1}$, otherwise we set $\mu^{+} \coloneqq + \infty$. Thus, we always have $\mu^{-} < \mu_{i_c} = 0 < \mu^+$.
		
		\begin{theorem}[Center manifold theorem]\label{center_manifold}
			Let the same assumptions as in Theorem \ref{stable_manifold} be satisfied. Then we find, for every $ 0 < \upsilon < \max\{\mu^+,\mu^-\} $, a family of immersed submanifolds $\mathcal{M}^{\upsilon,c}_{loc}(\omega)$ of $H_0$ and a set of full measure $\tilde{\Omega}$ of $\Omega$ with the following properties:
			\begin{enumerate}[(i)]
				\item There is a map $h_\omega^c \colon C_\omega\to\mathcal{M}^{\upsilon,c}_{loc}(\omega)$ with $h^c_{\omega}(0) = Z_{\omega}$ that is a homeomorphism, Lipschitz continuous and differentiable at zero. In other words, $\mathcal{M}^{\upsilon,c}_{loc}(\omega)$ is a topological Banach manifold modeled on $C_\omega$.
				
				\item $\mathcal{M}^{\upsilon,c}_{loc}(\omega)$  is $\tilde{\phi}$-locally invariant in the following sense: There is a random variable $\rho \colon \Omega \to (0,\infty)$ such that for
				\begin{align*}
					D^n_{\omega} \coloneqq \bigcap_{0 \leq j \leq n-1} \{Z_{\omega} + \xi \, \mid\, \|\tilde{\phi}_{\omega}^{j t_0}(Z_{\omega} + \xi) - Z_{\theta_{j t_0} \omega} \| \leq \rho(\theta_{(j+1)t_0} \omega) \},
				\end{align*}
				we have
				\begin{align*}
					\tilde{\phi}_{\omega}^{n t_0}(\mathcal{M}^{\upsilon,c}_{loc}(\omega) \cap D^n_{\omega}) \subset \mathcal{M}^{\upsilon,c}_{loc}(\omega)
				\end{align*}
				for every $n \geq 1$.
			\end{enumerate}
		\end{theorem}
		
		\begin{proof}
			Using the bounds that we obtained in the proof of Theorem \ref{stable_manifold}, the statements follow directly from \cite[Theorem 3.14 and Corollary 3.15]{Sebastian-maz23}.
		\end{proof}
		
		\begin{remark}
			The map $h_\omega^c$ and the immersed submanifolds $\mathcal{M}^{\upsilon,c}_{loc}(\omega)$ are given explicitly in \cite[Theorem 2.14]{Sebastian-maz23}.
		\end{remark}
		
		\subsection{Existence of stationary points.}
		The crucial condition in our invariant manifold theorems is the existence of a stationary point. Since we consider an SPDE with additive noise in \eqref{evolution-pb}, it is clear that a deterministic stationary point does not exist. In this section, we will prove that in the case of a bounded $G$, the existence of random stationary points can be deduced. This is, of course, a rather strong assumption, but our aim here is merely to illustrate that examples do exist where our theorems apply. We conjecture that the tools presented in \cite{FGS17a, FGS17b} will be useful when considering more general cases.

		\begin{proposition}\label{Example-G-bounded}
			Assume $\omega_{A}<0$ and that
			\begin{align*}
				\Vert G\Vert:=\sup_{\xi\in H_{0}}\Vert G(\xi)\Vert<\infty.
			\end{align*}
			Then there exists $M_{G}<\infty$ such that if $\Vert  G\Vert_{H}<M_{G}$, we can find a unique continuous map $V:\mathbb{R}\rightarrow H$ such that
			\begin{align}\label{nonlinear_stationary_point}
				\forall t\in\mathbb{R}: \ \ V(t) = \int_{-\infty}^{t}(-A_0)^{\beta}T_{0}(t-s)(-A)^{-\beta}G(V(s)+Y_{\theta_{s}\omega})\, \mathrm{d}s.
			\end{align}
			Furthermore, 
			\begin{align*}
				Z_{\omega} \coloneqq V(0) + Y_{\omega}
			\end{align*}
			is a stationary point for $\tilde{\phi}$ that satisfies the moment condition \eqref{eqn:moment_int_stat}.
		\end{proposition}

		\begin{proof}
			We start proving that \eqref{nonlinear_stationary_point} has a unique solution. Set
			\begin{align*}
				\mathcal{F}=\left\lbrace V\in C(\mathbb{R},H_0):\ \ \ \sup_{t\in\mathbb{R}}\Vert V(t)\Vert<\infty\right\rbrace.
			\end{align*}
			It is obvious that $\mathcal{F}$ is a Banach space.
			Define
			\begin{align}\label{fixed_point_stationary}
				\begin{split}
					&P:\mathcal{F}\rightarrow\mathcal{F},\\
					& P(V)(t):=\int_{-\infty}^{t} (-A_0)^{\beta}T_{0}(t-s)(-A)^{-\beta}G(V(s)+Y_{\theta_s\omega})\mathrm{d}s.
				\end{split}
			\end{align}
			We first claim that this map is well-defined. By assumption \eqref{fractional_growth},
			\begin{align}\label{well_define}	
				\begin{split}
					&\sup_{t\in\mathbb{R}}\left\Vert\int_{-\infty}^{t} (-A_0)^{\beta}T_{0}(t-s)(-A)^{-\beta}G(V(s)+Y_{\theta_s\omega})\mathrm{d}s\right\Vert\\&\quad\leq \sup_{t\in\mathbb{R}}\int_{-\infty}^{t}\left\Vert (-A_0)^{\beta}T_{0}(t-s)(-A)^{-\beta}G(V(s)+Y_{\theta_s\omega})\right\Vert \mathrm{d}s\leq M_{\beta}\Vert G\Vert\Vert (-A)^{-\beta}\Vert \int_{-\infty}^{0} (-s)^{-\beta}\exp(\omega_{A}(-s))\mathrm{d}s\\
					&\qquad<\infty. 
				\end{split}
			\end{align}
			Also by definition, $P(V)\in C(\mathbb{R},H_0)$, and the map is indeed well-defined. Assume $V_1,V_2\in \mathcal{F}$, then similar to \eqref{well_define}
			\begin{align*}
				&\sup_{t\in\mathbb{R}}\left\Vert P(V_2)(t)-P(V_1)(t)\right\Vert\leq\sup_{t\in\mathbb{R}}\int_{-\infty}^{t} \left\Vert(-A_0)^{\beta}T_{0}(t-s)(-A)^{-\beta}\big(G(V_2(s)+Y_{\theta_s\omega})-G(V_1(s)+Y_{\theta_s\omega})\big)\right\Vert\mathrm{d}s\\&\quad\leq \big(M_{\beta}\Vert G\Vert\Vert (-A)^{-\beta}\Vert \int_{-\infty}^{0} (-s)^{-\beta}\exp(\omega_{A}(-s))\mathrm{d}s\big)\times\sup_{t\in\mathbb{R}}\Vert V_{2}(t)-V_{1}(t)\Vert.
			\end{align*}
			Therefore, if 
			\begin{align*}
				\Vert G\Vert<\frac{1}{M_{\beta}\Vert (-A)^{-\beta}\Vert \int_{-\infty}^{0} (-s)^{-\beta}\exp(\omega_{A}(-s))\mathrm{d}s},
			\end{align*}
			running a standard fixed-point argument on \eqref{fixed_point_stationary}, we can find a unique $V\in \mathcal{F}$ such that
			\begin{align*}
				\forall t\in\mathcal{F}: \ \  V(t):=\int_{-\infty}^{t} (-A_0)^{\beta}T_{0}(t-s)(-A)^{-\beta}G(V(s)+Y_{\theta_s\omega})\mathrm{d}s.
			\end{align*}
			Clearly, $V(t)$ is random and we write $V_{\omega}(t)$ to make this visible. We claim that
			\begin{align}\label{eqn:V_stat}
				V_{\theta_t \omega}(s) = V_{\omega}(t + s).
			\end{align}
			Indeed, since
			\begin{align*}
				V_{\omega}(t + s) &= \int_{-\infty}^{t+s} (-A_0)^{\beta}T_{0}(t+s-\tau)(-A)^{-\beta}G(V_{\omega}(\tau)+Y_{\theta_{\tau}\omega}) \, \mathrm{d} \tau \\
				&= \int_{-\infty}^{s} (-A_0)^{\beta}T_{0}(s-\tau)(-A)^{-\beta}G(V_{\omega}(\tau + t) + Y_{\theta_{\tau + t}\omega}) \, \mathrm{d} \tau \\
				&= \int_{-\infty}^{s} (-A_0)^{\beta}T_{0}(s-\tau)(-A)^{-\beta}G(V_{\omega}(\tau + t) + Y_{\theta_{\tau} \theta_t \omega}) \, \mathrm{d} \tau
			\end{align*}
			and, by definition,
			\begin{align*}
				V_{\theta_t \omega}(s) = \int_{-\infty}^{s} (-A_0)^{\beta} T_{0}(s-\tau)(-A)^{-\beta} G(V_{\theta_t \omega} (\tau)+Y_{\theta_{\tau}\theta_t \omega})\, \mathrm{d} \tau,
			\end{align*}
			we see that $s \mapsto V_{\omega}(t + s)$ and $s \mapsto V_{\theta_t \omega}(s)$ solve the same equation, and \eqref{eqn:V_stat} follows by our uniqueness result. Recall that the solution $Z$ to  \eqref{evolution-pb} is given by 
			\begin{align*}
				Z_{\omega}(t) = V_{\omega}(t) + Y_{\omega}(t).
			\end{align*}
			With initial condition $Z_{\omega} \coloneqq Z_{\omega}(0) =  V_{\omega}(0) + Y_{\omega}(0)$, Lemma \ref{lemma:Y_stationary} and \eqref{eqn:V_stat} imply that
			\begin{align*}
				Z_{\theta_t \omega}(s) = Z_{\omega}(t + s)
			\end{align*}
			on a universal set of full measure. Therefore, if $t \geq 0$,
			\begin{align*}
				\tilde{\phi}^t_{\omega}(Z_\omega) = Z_{\omega}(t) = Z_{\theta_t \omega}(0) = Z_{\theta_t \omega}
			\end{align*}
			and we proved that $Z_{\omega}$ is indeed a stationary point. It remains to argue why $Z_{\omega}$ satisfies the moment condition \eqref{eqn:moment_int_stat}. To see this, it suffices to note that $Y_{\omega}(0)$ satisfies \eqref{eqn:moment_int_stat} since it is a stochastic integral with deterministic integrand and $V_{\omega}(0)$ has finite $p$-th moments for every $p \geq 1$ because $G$ is globally bounded. Thus, \eqref{eqn:moment_int_stat} follows for $Z_{\omega}$ by the triangle inequality.
			
			%         Let $\beta>1-\frac{1}{p^*}$ and $q^{*}\in \left[1,\frac{1}{1-\frac{1}{p^*}}\right)$, according to \eqref{regularity-integrated-sg} we have
			% \begin{align*}
				% \|\int_{-\infty}^{0}(-A_0)^{\beta}T_{0}(-s)(-A)^{-\beta}G(V(s)+Y_{\theta_{s}\omega})\mathrm{d}s\|&\leq\int_{-\infty}^{0}\|(-A_0)^{\beta}T_{0}(-s)(-A)^{-\beta}G(V(s)+Y_{\theta_{s}\omega})\|\mathrm{d}s\\&\leq
				% (M_{\beta}\Vert G\Vert\Vert (-A)^{-\beta}\Vert)^{q^*} \int_{-\infty}^{0} (-s)^{-\beta q^*}\exp(\omega_{A}(-s))\mathrm{d}s\\
				% &\leq
				% (M_{\beta}\Vert G\Vert\Vert (-A)^{-\beta}\Vert)^{q^*} \int^{+\infty}_{0} (s)^{-\beta q^*}\exp(\omega_{A}s)\mathrm{d}s
				% \end{align*}
			% Since $\omega_A <0$, then the last integral converges, consequently 
			% \begin{align*}
				% \mathbb{E}\|V(0)\|^p\leq
				% (M_{\beta}\Vert G\Vert\Vert (-A)^{-\beta}\Vert)^{pq^*} \left(\int^{+\infty}_{0} (s)^{-\beta q^*}\exp(\omega_{A}s)\mathrm{d}s\right)^p<\infty.
				% \end{align*}

			% Finally if we set $Z(t):=V(t)+Y_{\theta_{t}\omega}$, then $Z$ the stationary point satisfying in \eqref{nonlinear_stationary_point}. The claim follows from the fact that $\tilde{\phi}^{t}_{\omega}(Z_{\omega})=Z_{\theta_t\omega}$, where $Z_{\omega}$ is the initial data.
		\end{proof}
		
		\begin{comment}
			\begin{definition}
				For $0<\mu<-\mu_{i_0}$, we set
				\begin{align*}
					BC^{\mu}:=\lbrace V\in C[0\,\infty):\sup_{s\geq 0}\exp(-\mu s)\Vert V(s)\Vert<\infty\rbrace
				\end{align*}
			\end{definition} 
			\begin{lemma}
				Following map is well defined
			\end{lemma}
			\begin{align*}
				&\mathcal{F}_{\omega}: S_{\omega}\times BC^{\mu}\longrightarrow BC^{\mu},\\ &\quad
				\mathcal{F}_{\omega}(\xi,V):\psi^{t}_{\omega}(\xi)+\int_{0}^{t}(-A_0)^{\beta}\psi^{t-s}_{\theta_s\omega}\Pi_{S_{\theta_s\omega}}(-A)^{-\beta}\big[G(V(s)+Z_{\theta_s\omega})-G(Z_{\theta_{s}\omega})-D_{Z_{\theta_s\omega}}G(V(s))\big]\mathrm{d}s\\&\qquad-\int_{t}^{\infty}(-A_0)^{\beta}\psi^{t-s}_{\theta_s\omega}\Pi_{U_{\theta_s\omega}}(-A)^{-\beta}\big[G(V(s)+Z_{\theta_s\omega})-G(Z_{\theta_{s}\omega})-D_{Z_{\theta_s\omega}}G(V(s))\big]\mathrm{d}s
			\end{align*}
			\begin{proof}
				details later
			\end{proof}
		\end{comment}
		\section{Examples: Heat Equation with white noise at the boundary}\label{Application-section}
		Let ${O}\subset\mathbb{R}^{n}$ be an open bounded smooth domain with boundary $\partial{O}=:\Gamma$ and $\mu>0$. Consider the following stochastic SPDE with white noise at the boundary
		\begin{align}\label{Neumann-Example}
			\begin{cases}
				\frac{\partial u(t,x)}{\partial t}=\frac{\partial^{2}u(t,x)}{\partial x^{2}}+\mu u(t,x)+g(u(t,x)),&\quad t>0,\quad x\in{O}\\
				\partial_\nu u(t,x)=\dot{W}_{0}(t,x),&\quad t>0,\quad x\in\Gamma\\
				u(0,\cdot)=\xi.
			\end{cases}
		\end{align}
		Here $\nu=\nu(x)$ is the unit normal vector of $\Gamma$ pointing towards the exterior of ${O}$, $\dot{W}_{0}(t)$, $t\geq0$ is a $Q$-Wiener process on the probability space $(\Omega,\mathcal{F},\mathbb{P})$, $g:L^{2}({O})\to L^{2}({O})$ is a Lipschitz continuous function assumed to be of class $C^1$ and globally bounded and the initial state $\xi\in L^2(O)$. In order to deal with boundary conditions, we define the product space
		\begin{equation*}
			{H}:=L^{2}(\Gamma)\times L^{2}({O}),
		\end{equation*}
		as well as we select the following operator
		\begin{align*}
			D({A}):=\left\{{0}\right\}\times H^{2}({O}),\quad {A}\begin{pmatrix}0\\\phi\end{pmatrix}:=\begin{pmatrix}\partial_\nu\phi\\\varphi^{\prime \prime}\end{pmatrix}.
		\end{align*}
		One may observe that
		\begin{equation*}
			{H}_{0}=:	\overline{D({A})}=\left\{{0}\right\}\times L^{2}({O})\neq H,
		\end{equation*}
		So that the operator $A$ is non-densily defined. Moreover, the part of $A$ in $H_0$
		
		$$
		A_0\left(\begin{array}{l}
			0 \\
			\varphi
		\end{array}\right)=\left(\begin{array}{c}
			
			0 \\
			\varphi^{\prime \prime}
		\end{array}\right) \text { with } D\left(A_0\right)=\{0\}\times\left\{\varphi \in H^{2}(O): \partial_\nu\varphi=0\right\}
		$$
		is an infinitesimal generator of a strongly continuous analytic semigroup $\{T_0(t)\}_{t\geq 0}$ on $H$.\\
		In the context of parabolic equations, it is well known that elliptic operators, due to the boundary conditions, are not generally Hille-Yosida operators but almost sectorial operators (see \cite{Daprato-66}). In particular, according to \cite[Theorem 2.1]{Agra-Denk-97}, the resolvent of $A$ satisfies the estimate
		\begin{align*}
			&0<\liminf _{\lambda \rightarrow \infty} \lambda^{\frac{1}{p^*}}\left\|(\lambda I-A)^{-1}\right\|_{\mathcal{L}(H)} \leq \limsup _{\lambda \rightarrow \infty} \lambda^{\frac{1}{p^*}}\left\|(\lambda I-A)^{-1}\right\|_{\mathcal{L}(H)}<\infty, \text { where } p^*=\frac{4}{3}.
		\end{align*}
		Thus, the operator $A$ is $p^*$-sectorial operator, which ensures that the assumption \ref{Assumption1.1} are satisfied.
		
		Now by taking 
		\begin{align*}
			G(X(t)):=\begin{pmatrix}0\\g(u(t,\cdot))\end{pmatrix},\quad W(t):=\begin{pmatrix}W_{0}(t,\cdot)\\0\end{pmatrix},\quad t>0,
		\end{align*}
		and identifying the function $u(t,\cdot)$ with $V(t)=\begin{pmatrix}0\\u(t,\cdot)\end{pmatrix}$, the boundary problem \eqref{Neumann-Example} takes the following form  
		\begin{align}\label{example-pb}
			\begin{cases}
				\d V(t)=(A+\mu I)V(t) \, \mathrm{d}t+G(V(t)) \, \mathrm{d}t +  \mathrm{d}W(t), \quad t\geq0,\cr V(0)=V_0\in H_0,
			\end{cases}
		\end{align}
		where $I$ stands for the identity operator on $H$. Furthermore, the spectrum of $A+\mu I$ satisfies
		\begin{align*}
			\sigma(A+\mu I)=\sigma(A_0+\mu I)=\{\mu-\pi_n:n\in\mathbb{N}\},
		\end{align*}
		where $\pi_n>0$ are only isolates points, namely the spectrum $\sigma(A_0)=\{-\pi_n:n\in\mathbb{N}\}$.\\
		Consequently $A_{0,\mu}:=A_0+\mu I$ generates a compact, analytic $C_0$-semigroup $(T_{0,\mu}(t))_{t \geq 0}$ on $H_0$, and
		
		$$
		\|T_{0,\mu}(t)\|_{\mathcal{L}(H_0)} \leq e^{-\mu t}, \text { for } t \geq 0.
		$$
		
		It follows from Theorem \ref{Existence-solution} and Corollary \ref{dynamical-sys2} that the boundary parabolic problem \eqref{Neumann-Example} has a unique integrated solution which generates a random dynamical system $\tilde{\phi}:\mathbb{R}^2\times \Omega\times H_0\to H_0$. As $g$ is assumed to be of class $C^1$, then according to Proposition \ref{Fréchet}, $\tilde{\phi}$ is Fréchet differentiable and the derivative satisfies for any $\xi,\eta\in H_0$ and $t>0$ and $\lambda\in\rho(A)$
		\begin{equation}\label{Fréchet_differentiable2}
			D_{\xi}\tilde{\phi}^t_{\omega}[\eta]=T_{0,\mu}(t)\eta+\lim_{\lambda\to\infty}\int_0^tT_{0,\mu}(t-s)\lambda(\lambda-A)^{-1}DG_{\tilde{\phi}^s_{\omega}(\xi)}(D_{\xi}\tilde{\phi}^s_{\omega}[\eta])\mathrm{d}s.
		\end{equation}
		
		The compactness of \eqref{Fréchet_differentiable2} comes from the fact that the semigroup $(T_{0,\mu}(t))_{t\geq0}$ and the resolvent operator $(\lambda-A)^{-1}$ are compact. It follows from Proposition \ref{compactness-result} that the linearized solution $D_{\xi}\tilde{\phi}^t_{\omega}$ given by \eqref{Fréchet_differentiable2} is indeed compact. To ensure the existence of stationary solution for $\tilde{\phi}$, we need to verify the moment condition \ref{eqn:moment_int_stat}. 
		%First, as the semigroup $(T_{0}(t))_{t\geq0}$ is compact and the associated spectre is discrete, then it can be represented as a series expansion in terms of the eigenfunctions $\phi_k$ and their corresponding eigenvalues $\lambda_k$ as follows
		% \begin{align*}
			%     T_{0}(t)\xi=\sum_{k=0}^{\infty} e^{-\lambda_k t} <\xi,\phi_k>\phi_k.
			%  \end{align*}
		% It follows that 
		%\begin{align*}
		%\left\|\frac{d S_A(t) \xi}{d t}\right\| & =\left\|\left(-A_0\right)^\beta T_{A_0}(t)(-A)^{-\beta} \xi\right\|,\\
		%& \leq \left\|(-A)^{-\beta}\right\| \| \sum_{k=0}^{\infty} \lambda_k^\beta e^{-\lambda_k t}\langle\xi, \phi_k\rangle\phi_k \|,\\
		%& \leq M_\beta \left\|(-A)^{-\beta}\right\|\|x\|,
		%\end{align*}
		As $g$ is assumed to be globally bounded, then from Proposition \ref{Example-G-bounded} the stationary point exists and it's given by \eqref{nonlinear_stationary_point}. Now, the existence of invariant manifolds (stable, unstable, and center) around the stationary point follows immediately from Theorem \ref{stable_manifold}, Thoerem \ref{unstable_manifold}, Theorem \ref{center_manifold} and Proposition \ref{Example-G-bounded}.
		
		\begin{remark}
			The preceding results can be extended to the cylindrical case, i.e, when $Q=I$, specifically for the one-dimensional space ($n=1$). In fact, we know that for any $\beta\in(\frac{1}{4},1)$ and $t>0$ the derivative of the integrated semigroup is given by
			\begin{align*}
				\frac{d}{dt}{S}_{{A}}(t)x&=(-{A}_0)^{\beta}{T}_{0}(t)(-A)^{-\beta}x\\
				&=(\lambda_k)^{\beta}e^{-\lambda_k t}(\lambda_k)^{-\beta}x,
			\end{align*}
			where $(\lambda_k)_{k\geq1}$ denoted the eigenvalues of $-{A}_0$. Then for any $\tau>0$
			\begin{align*}
				\displaystyle\int_{0}^{\tau}\left\|\frac{d{S}_{{A}}(r)}{dr}\right\|_{{L}_2({H})}^{2}dr&=\displaystyle\int_{0}^{\tau} Tr ((-{A}_0)^{2\beta}e^{2r A_0}(-{A})^{-2\beta}) dr\\
				&=\sum_{k=1}^{\infty}\frac{1-e^{-2\lambda_k\tau}}{2\lambda_k}.
			\end{align*}
			Since $\lambda_k\sim ck^{\frac{2}{n}}$, the above series convergent if and only if $n=1$. 
		\end{remark}

		\appendix
		\section{}\label{appendix}
		The following Gronwall inequality is crucial for our main results.
		\begin{lemma}\label{powered_Gronwall-type_inequality}
			Let $0< \beta<1$ and $T>0$ and $\kappa: [0,T]\rightarrow (0,\infty)$, be a nondecreasing continuous function. Also, $g:[0,T]\rightarrow (0,\infty)$ is locally integrable. Assume $u:[0,T]\rightarrow(0,\infty)$ is a locally integrable function such that
			\begin{align*}
				t\in [0,T]:\ \ \  u(t)\leq g(t)+\kappa(t)\int_{0}^{t}(t-s)^{-\beta}u(s)\mathrm{d}s.
			\end{align*}
			Then
			\begin{align*}
				t\in [0,T]:\ \ \  u(t)\leq \int_{0}^{t}\sum_{n\geq 1}\frac{[\kappa(t)\Gamma(1-\beta)]^n}{\Gamma(n(1-\beta))}(t-s)^{n(1-\beta)-1}g(s)\mathrm{d}s.
			\end{align*}
		\end{lemma}
		\begin{proof}
			Cf.\cite[Theorem 1]{YE20071075}.
		\end{proof}
		\begin{remark}
			Recall that the Gamma function $\Gamma(n+1)=n!$ and for any positive real number $x$ $\Gamma(x+1)=x\Gamma(x)$. Moreover, for some $\mu\in(1,2)$ this function on $(0,\mu]$ is decreasing and on $(\mu,\infty)$ is increasing.
		\end{remark}
		\begin{lemma}\label{sum_bound}
			Assume $\alpha\geq 1$, then
			\begin{align}\label{sum_bound_1}
				\sum_{n\geq 1}\frac{\alpha^{n-1}}{\Gamma(n(1-\beta))}\leq\frac{\beta+1}{1-\beta}\alpha^{\frac{1+\beta}{1-\beta}}+\frac{\alpha^{\frac{2}{1-\beta}}\exp(\alpha^{\frac{1}{1-\beta}})}{1-\beta}:=R(\alpha).
			\end{align}
			
			\begin{proof}
				Let $\alpha\geq 1$ and $\beta>1-\frac{1}{p^*}$, we obtain
				\begin{align*}
					\sum_{n\geq 1}\frac{\alpha^{n-1}}{\Gamma(n(1-\beta))}&=\sum_{1\leq n <\frac{\mu}{1-\beta}}\frac{\alpha^{n-1}}{\Gamma(n(1-\beta))}+\sum_{ n \geq\frac{\mu}{1-\beta}}\frac{\alpha^{n-1}}{\Gamma(n(1-\beta))}\\&\leq \frac{\beta+\mu-1}{1-\beta}\alpha^{\frac{\beta+\mu-1}{1-\beta}}+\sum_{ n \geq\frac{\mu+\beta-1}{1-\beta}}\frac{(\alpha^{\frac{1}{1-\beta}})^{n(1-\beta)}}{\Gamma((n+1)(1-\beta))}\\&\leq \frac{\beta+1}{1-\beta}\alpha^{\frac{1+\beta}{1-\beta}}+\sum_{m\geq 1}\sum_{\frac{m+\beta-1}{1-\beta}\leq n<\frac{m+\beta}{1-\beta}}\frac{\alpha^{\frac{m+1}{1-\beta}}}{(m-1)!}\\&\leq \frac{\beta+1}{1-\beta}\alpha^{\frac{1+\beta}{1-\beta}}+\frac{\alpha^{\frac{2}{1-\beta}}\exp(\alpha^{\frac{1}{1-\beta}})}{1-\beta}.
				\end{align*}
			\end{proof}
		\end{lemma}
		\bibliographystyle{plain}
		\bibliography{bibliography}
		
		% \begin{thebibliography}{99}
			% \bibitem{Ducrot-Magal-Prevost}
			% Ducrot, A., Magal, P., \& Prevost, K. (2010). \emph{Integrated semigroups and parabolic equations}. Part I: linear perburbation of almost sectorial operators. Journal of Evolution Equations, {\bf 10} (2), 263-291.
			% \bibitem{Ducrot-Lahbiri} Ducrot, A., \& Lahbiri, F. Z. (2023). \emph{Abstract Parabolic Equations with boundary white noise: an integrated semigroup approach}. arXiv preprint arXiv:2307.12294.
			% \bibitem{Magal-Ruan}
			%  Magal, P., \& Ruan, S. (2018). \emph{Semilinear Cauchy Problems with Non-dense Domain}. In Theory and Applications of Abstract Semilinear Cauchy Problems (pp. 217-248). Springer, Cham.
			%  \end{thebibliography}
	\end{document}